\documentclass[12pt]{article}
\usepackage{e-jc}
\usepackage{cmap}
\usepackage{mathtools}
\usepackage{amsmath}
\usepackage{amsfonts}
\usepackage{amssymb}
\usepackage{amsthm}
\usepackage{enumerate}
\usepackage[english]{babel}
\usepackage[T1]{fontenc}
\usepackage{graphicx,epsfig,color}
\usepackage[colorlinks=true,citecolor=black,linkcolor=black,urlcolor=blue]{hyperref}
\usepackage{cleveref}
\usepackage{float}
\usepackage[mathscr]{euscript}
\usepackage{ifpdf}
\usepackage{caption}
\usepackage{cite}
\usepackage{tikz-cd}
\usepackage[toc,page]{appendix}
\usepackage{listings}

\newcommand{\bs}{\ensuremath{\backslash}}

\newcommand{\arxiv}[1]{\href{http://arxiv.org/abs/#1}{\texttt{arXiv:#1}}}

\theoremstyle{plain}
\newtheorem{theorem}{Theorem}
\newtheorem{lemma}[theorem]{Lemma}
\newtheorem{proposition}[theorem]{Proposition}

\newtheorem{algorithm}[theorem]{Algorithm}

\newtheorem{construction}[theorem]{Construction}

\theoremstyle{definition}
\newtheorem{Definition}[theorem]{Definition}
\newtheorem{Remark}[theorem]{Remark}
\newtheorem{Notation}[theorem]{Notation}
\newtheorem{Example}[theorem]{Example}

    \newtheoremstyle{TheoremNum}
        {\topsep}{\topsep}              
        {\itshape}                      
        {}                              
        {\bfseries}                     
        {.}                             
        { }                             
        {\thmname{#1}\thmnote{ \bfseries #3}}
    \theoremstyle{TheoremNum}

\usepackage{blindtext, color}
\newcommand{\PC}{\mathbb{P}^1(\mathbb{C})}

\title{Pruned double Hurwitz numbers}
\author{Marvin Anas Hahn\\
\small Mathematisches Institut\\[-0.8ex]
\small Universit\"at T\"ubingen\\[-0.8ex] 
\small Auf der Morgenstelle 10, 72076 T\"ubingen, Germany\\
\small\tt marvin-anas.hahn@uni-tuebingen.de\\
}

\date{\dateline{Jan 1, 2012}{Jan 2, 2012}}

\begin{document}
\maketitle

\begin{abstract}
 Hurwitz numbers count ramified genus $g$, degree $d$ coverings of the projective line with fixed branch locus and fixed ramification data. Double Hurwitz numbers count such covers, where we fix two special profiles over $0$ and $\infty$ and only simple ramification else. These objects feature interesting structural behaviour and connections to geometry. In this paper, we introduce the notion of pruned double Hurwitz numbers, generalizing the notion of pruned simple Hurwitz numbers in \cite{DN13}. We show that pruned double Hurwitz numbers, similar to usual double Hurwitz numbers, satisfy a cut-and-join recursion and are piecewise polynomial with respect to the entries of the two special ramification profiles. Furthermore double Hurwitz numbers can be computed from pruned double Hurwitz numbers. To sum up, it can be said that pruned double Hurwitz numbers count a relevant subset of covers, leading to considerably smaller numbers and computations, but still featuring the important properties we can observe for double Hurwitz numbers.
 
\end{abstract}

\section{Introduction}
\label{chapter:Introduction}
\pagenumbering{arabic}
Hurwitz numbers are important enumerative objects connecting numerous areas of mathematics, such as algebraic geometry, algebraic topology, operator theory, representation theory of the symmetric group and combinatorics. Historically, these objects were introduced by Adolf Hurwitz in \cite{Hurwitz} to study the moduli space $\mathcal{M}_g$ of curves of genus $g$.\\
There are various equivalent definitions of Hurwitz numbers and several different settings, among which the most well-studied one is the case of simple Hurwitz numbers, which we denote by $\mathcal{H}_g(\mu)$. To be more precise, simple Hurwitz numbers count genus $g$ coverings of $\PC$ with fixed ramification $\mu$ over $0$ and simple ramification over $r$ further fixed branch points, where the number $r$ is given by the Riemann-Hurwitz formula. The theory around these objects is well developed and a lot is known about their structure. Each degree $d$ cover $f:Y\to X$ with branch locus $B$ induces a \textit{monodromy representation}, i.e. a map $\phi:\pi_1(X\backslash B)\to \mathbb{S}_d$. Starting from these monodromy representations and applying Riemann's existence theorem one can show that there is an equivalent definition in terms of factorizations of permutations (see chapter 7.2 in \cite{CM}). Moreover, simple Hurwitz numbers satisfy a cut-and-join recursion which is inherent in the combinatorial structure of these factorizations. Another well-known result is the fact that --- up to a combinatorial factor --- $\mathcal{H}_g(\mu)$ behaves polynomially in the entries of $\mu$ for fixed genus $g$ and fixed length of $\mu$. Recently, there has been an increased interest in Hurwitz theory due to connections to Gromov-Witten theory, remarkably through the celebrated ELSV formula \cite{ekedahl2001hurwitz} which relates Hurwitz numbers to intersection products in the moduli space of curves. This formula initiated a rich interplay between those areas. The polynomiality result for simple Hurwitz numbers is a consequence of the ELSV formula. Via the ELSV formula a new proof of Witten's conjecture was given in \cite{OPsecondversion} using Hurwitz theory. Moreover, simple Hurwitz numbers satisfy the Chekhov-Eynard-Orantin topological recursion, a theory motivated by mathematical physics with numerous applications in geometry (see e.g. \cite{chekhov2006hermitian}, \cite{BouMa}, \cite{EMS}, \cite{MR2519749}).\vspace{\baselineskip}\\
A further case which has been of great interest in recent years is the one of double Hurwitz numbers, which we denote by $\mathcal{H}_g(\mu,\nu)$. Here we allow two special ramification profiles, that is in addition to allowing arbitrary ramification $\mu$ over $0$, we allow arbitrary ramification $\nu$ over $\infty$. Obviously, for $\nu=(1,\dots,1)$ this yields the definition for simple Hurwitz numbers given above. While there are still a lot of open questions, much is known about these objects as well and they admit many results, which are similar to those about simple Hurwitz numbers. Among those is a cut-and-join recursion for double Hurwitz numbers and a definition in terms of factorizations in the symmetric group. In \cite{GJV} it was proved that $\mathcal{H}_g(\mu,\nu)$ behaves piecewise polynomially in the entries of $\mu$ and $\nu$. More than that, wall-crossing formulas in genus $0$ were given in \cite{SSV} and in all genera in \cite{CJMa} and \cite{johnson2015double}.\\
Among the open problems for double Hurwitz numbers is the question, if there is an ELSV-type formula for them \cite{cavalieri2016hurwitz}. Some progress has been made in \cite{GJV}, where such a formula is given for genera $0$ and $1$. Furthermore, it is not known, whether double Hurwitz numbers satisfy an Eynard-Orantin topological recursion.\vspace{\baselineskip}\\
In \cite{DN13} the notion of pruned simple Hurwitz numbers was introduced. The main idea behind this notion is, that it is sufficient to consider a non-trivial subset of ramified covers that contribute to the simple Hurwitz number that still carries all the information and that this subset may be described purely combinatorially in terms of certain graphs on surfaces. These graphs were introduced as branching graphs in \cite{OPsecondversion}. There are various names in the literature for these and similar graphs, such as ribbon graphs, dessin d'enfants, Hurwitz galaxies, maps in surfaces, graphs in surfaces. The pruned simple Hurwitz number, which we denote by $\mathcal{PH}_g(\mu)$ is a count over this subset. It was established in \cite{DN13}, that simple Hurwitz numbers and pruned simple Hurwitz numbers are equivalent in the sense, that simple Hurwitz numbers may be computed as a weighted sum over certain pruned simple Hurwitz numbers of the same genus. Moreover, these new objects still carry a lot of information of the standard case, such as the fact that $\mathcal{PH}_g(\mu)$ behaves polynomially in the entries of $\mu$. Pruned simple Hurwitz numbers are defined in terms of graphs on surfaces, however there is a definition in terms of factorizations of permutations, as well. Moreover, they admit a cut-and-join recursion similar to the one for simple Hurwitz numbers. Using these results and the ELSV formula, another proof for Witten's Conjecture was given in \cite{DN13}. Furthermore, it was proved, that pruned simple Hurwitz numbers admit an Eynard-Orantin topological recursion.

To sum up, it can be said that simple pruned Hurwitz numbers count a relevant subset of covers, leading to considerably smaller numbers and computations, but still featuring the important properties we can observe for simple Hurwitz numbers.\vspace{\baselineskip}\\
The aim of this paper is to introduce the notion of pruned double Hurwitz numbers, generalizing the definition in \cite{DN13} and to investigate their structure. Our definition of pruned double Hurwitz numbers, which we denote by $\mathcal{PH}_g(\mu,\nu)$, is given in terms of branching graphs, as well. We prove three structural results about pruned double Hurwitz numbers:
\begin{theorem}
\label{1}
 Double Hurwitz numbers can be expressed in terms of pruned double Hurwitz numbers with smaller input data (i.e. smaller degree and ramification data, but the same genus).
\end{theorem}
 For a precise formulation see Theorem \ref{thm:mainresult}. \Cref{chapter:PrunedHurwitz} is devoted to the proof of this theorem.
\begin{theorem}
\label{2}
 Pruned double Hurwitz numbers satisfy a cut-and-join recursion.
\end{theorem}
 For a precise formulation see Theorem \ref{thm:prunedrecursion}, which is proved in \Cref{sub:rec}.
\begin{theorem}
\label{3}
 Pruned double Hurwitz numbers are piecewise polynomial with the same chamber structure as in the standard case.
\end{theorem}
For a precise formulation see Theorem \ref{thm:poly}, which is proved in the first half of \Cref{sub:poly}.

Moreover, we express pruned double Hurwitz numbers in terms of factorizations in the symmetric group. We begin this paper by recalling some basic facts about Hurwitz numbers and re-introducing branching graphs in a way suitable for our purposes in Section \ref{chapter:Preliminaries}. In Section \ref{sub:corr}, we introduce the notion of pruned double Hurwitz numbers and prove Theorem \ref{1}. We continue in Section \ref{sub:rec} by formulating and proving Theorem \ref{2}. In Section \ref{sub:poly}, we give a proof for Theorem \ref{3}. We note, that while our first two results are proven in a similar way as their corresponding results in \cite{DN13}, the method used for the polynomiality result is not feasible for pruned double Hurwitz numbers. In fact, our method is similar to the one used in \cite{GJV} to prove the piecewise polynomiality for double Hurwitz numbers. We finish this section by connecting the combinatorics of branching graphs to the combinatorics of symmetric groups and express pruned double Hurwitz numbers in the setting of factorizations of permutations. Building on these results, we developed and implemented an algorithm to compute pruned double Hurwitz numbers. An implementation of the algorithm in the computer algebra system \cite{GAP4} may be found in \href{https://sites.google.com/site/marvinanashahn/computer-algebra}{https://sites.google.com/site/marvinanashahn/computer-algebra}. Using this tool, we computed several non-trivial examples of Hurwitz numbers and pruned Hurwitz numbers. The computations agree with the predictions made by the formulas of Theorem \ref{thm:mainresult} and Theorem \ref{thm:prunedrecursion}.

\section{Preliminaries}
\label{chapter:Preliminaries}
In this section, we introduce some basic notions of graph theory and the theory of Hurwitz numbers. Detailed introductions to these topics can be found in \cite{RS}, \cite{RM} p.84-92 and the book \cite{CM}.

\subsection{Graphs}
We consider graphs with half edges $(V,E,E')$. Here $V$ is the set of vertices and the multiset $E\subset V\times V$ is the set of edges. The multiset $E'\subset V$ is the set of half-edges. A \textit{forest} is a graph without cycles and a \textit{tree} is a connected forest.

We note, that we define the valency $\mathrm{val}(v)$ to be the number of full-edges incident to $v$. By convention, we count loops twice.
\begin{figure}[H]
\begin{center}
\scalebox{0.7}{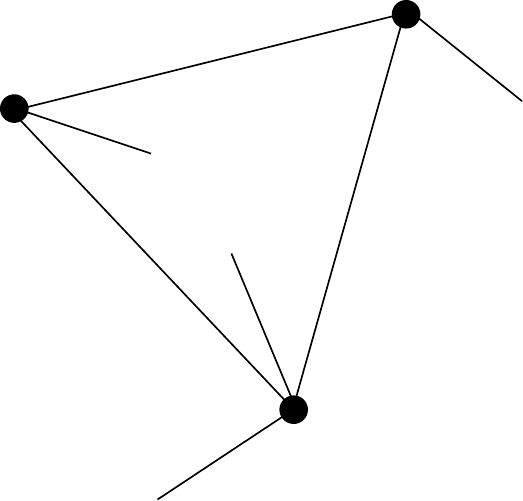}
\end{center}
\caption{A graph with half-edges.}
\end{figure}

\begin{figure}
 \begin{center}
  \scalebox{.4}{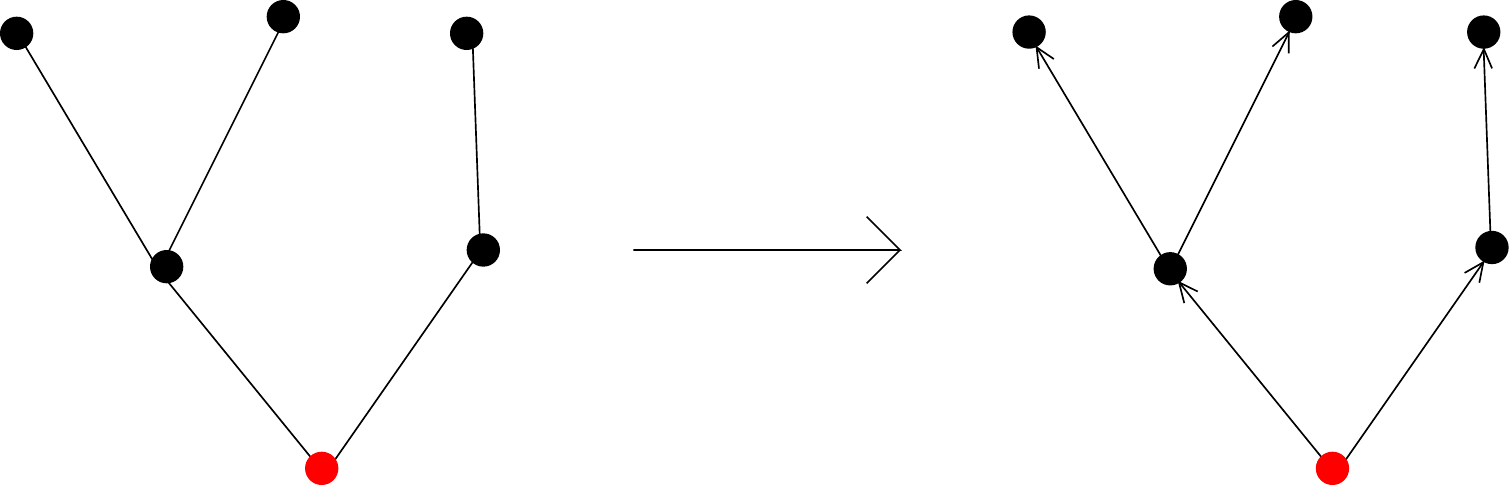}
 \end{center}
\caption{Rooted forest (root is marked red) with canonical orientation.}
\label{fig:root}
\end{figure}

Obviously, we may decompose each graph into its connected components. We call a forest \textit{rooted}, if each component contains a distinguished vertex, which we call the \textit{root-vertex}. Note that a rooted forest carries a canonical orientation in the way, that the edges of each connected component point away from the corresponding root-vertex (see e.g. Figure \ref{fig:root}).

\begin{Definition}
Let $v$ be a vertex in a rooted forest with the canonical orientation. We call the target vertex of an outgoing edge of $v$ a \textit{successor} of $v$.
\end{Definition}

\subsection{Hurwitz numbers}

\begin{Definition}
 Let $d$ be a positive integer, $\mu,\nu$ two ordered partitions of $d$ and let $g$ be a non-negative integer. Moreover, let $p_1,p_2,q_1,\dots,q_m$ be points in $\mathbb{P}^1(\mathbb{C})$, such that $m=2g-2+\ell(\mu)+\ell(\nu)$. We define a Hurwitz cover of type $(g,\mu,\nu)$ to be a branched covering $f:C\to\PC$, such that:
 \begin{enumerate}
  \item $C$ is a genus $g$ curve,
  \item $f$ is a degree $d$ map, that ramifies with profile $\mu$ over $p_1$, with profile $\nu$ over $p_2$ and $(2,1\dots,1)$ over $q_i$,
  \item $f$ is unramified everywhere else,
  \item the pre-images of $p_1$ and $p_2$ are labeled by $1,\dots,\ell(\mu)$ and $1,\dots,\ell(\nu)$ respectively, such that the point labeled $i$ in $f^{-1}(p_1)$ (respectively $f^{-1}(p_2)$) has ramification index $\mu_i$ (respectively $\nu_i)$.  
 \end{enumerate}
 We call a branch point with ramification profile $(2,1,\dots,1)$ a simple branch point and we call a ramification point with ramification index $2$ a simple ramification point.\\
 An isomorphism between two covers $f:C\to\mathbb{P}^1(\mathbb{C})$, $f':C'\to\mathbb{P}^1(\mathbb{C})$ is a homeomorphism $\pi:C\to C'$ respecting the labels, such that $f'\circ \pi=f$. We denote the automorphism group of a cover $f$ by $\mathrm{Aut}(f)$.\\ 
 Let $\mathbb{H}_g(\mu,\nu)$ be the set of all Hurwitz covers of type $(g,\mu,\nu)$. Then we define the double Hurwitz number \[\mathcal{H}_g(\mu,\nu)=\sum_{f\in\mathbb{H}_g(\mu,\nu)}\frac{1}{|\mathrm{Aut}(f)|}.\]                                                                                                                                                                                                                                                                                                
 Note that $\mathcal{H}_g(\mu,\nu)$ is a topological invariant, that is, it is independent of the locations of the points $p_1,p_2,q_1,\dots,q_m$ and of the complex structure of $C$.
 \end{Definition}
 
  By matching a cover with a monodromy representation, we may count ramified coverings of $\mathbb{P}^1(\mathbb{C})$ in terms of factorizations of permutations. For a permutation $\sigma$, denote by $\mathcal{C}(\sigma)$ the corresponding partition given by its decomposition in disjoint cycles.
\begin{theorem} \label{thm:symmetricgroup}Let $\mu$ and $\nu$ be ordered partitions of some positive integer $d$. Moreover, let $g$ be some non-negative integer. The following equation holds:
  \begin{align*}
   &\mathcal{H}_g(\mu,\nu)=\\&\frac{1}{d!}\left|
   \begin{cases}
   \begin{rcases}
   &\left(\sigma_1,\tau_1,\dots,\tau_m,\sigma_2\right)\textrm{, such that:}\\ 
   &\bullet\ \sigma_1,\ \sigma_2,\ \tau_i\in\mathcal{S}_d,\\
   &\bullet\ \sigma_2\cdot\tau_m\cdot\dots\cdot\tau_1\cdot\sigma_1=\mathrm{id},\\
   &\bullet\ \mathcal{C}(\sigma_1)=\mu,\ \mathcal{C}(\sigma_2)=\nu\textrm{ and }\mathcal{C}(\tau_i)=(2,1,\dots,1),\\
   &\bullet\ \textrm{the group generated by }\left(\sigma_1,\tau_1,\dots,\tau_m,\sigma_2\right)\textrm{ acts transitively on}\\ &\ \{1,\dots,d\},\\
   &\bullet\ \textrm{the disjoint cycles of }\sigma_1\textrm{ and }\sigma_2\textrm{ are labeled,}\textrm{ by }1,\dots,\ell(\mu)\textrm{ and }\\
   &\ 1,\dots,\ell(\nu)\textrm{ respectively},\\
   &\bullet \textrm{ the cycle of }\sigma_1\textrm{ (resp. }\sigma_2\textrm{) labeled }i\textrm{ has length }\mu_i\textrm{ (resp. }\nu_i\textrm{)}. 
   \end{rcases}
   \end{cases}
\right|.
  \end{align*}
\end{theorem}

\begin{proof}
 For a proof, see for example \cite{RM}.
\end{proof}

\subsection{Hurwitz galaxies and Branching graphs}
In this subsection, we explain a connection between covers contributing to $\mathcal{H}_g(\mu,\nu)$ and graphs on surfaces. We will define two notions of graphs on surfaces, that will turn out to be equivalent. We will start by defining branching graphs. We note that we will view full-edges as two half-edges glued together at their respective vertex-free ends.

 \begin{Definition}
\label{def:branchinggraph}
 Let $d$ be a positive integer, $\mu$ and $\nu$ be ordered partitions of $d$. We define a branching graph of type $(g,\mu,\nu)$ to be a graph $\varGamma$ embedded on an oriented surface $S$ of genus g, such that for $m=2g-2+\ell(\mu)+\ell(\nu)$:
 \begin{enumerate}[(i)]
 \item $S\backslash\varGamma$ is a disjoint union of open disks.
  \item There are $\ell(\mu)$ vertices, labeled $1,\dots,\ell(\mu)$, such that the vertex labeled $i$ is adjacent to $\mu_i\cdot m$ half-edges, labeled cyclically counterclockwise by $1,\dots,m$. We define the perimeter of the vertex labeled $i$ by $per(i)=\mu_i$.
  \item There are exactly $m$ full edges labeled by $1,\dots,m$.
  \item The $\ell(\nu)$ faces are labeled by $1,\dots,\ell(\nu)$ and the face labeled $i$ has perimeter $per(i)=\nu_i$, by which we mean, that each label occurs $\nu_i$ times inside the corresponding face, where we count full-edges adjacent to $i$ on both sides twice.
 \end{enumerate}
Note that we allow loops at the vertices. An isomorphism between two Hurwitz galaxies, is an orientation-preserving homeomorphism of their respective surfaces, which induces an isomorphism of graphs, that preserves vertex-, (half-)edge- and face-labels.
\end{Definition}

Now we will define a second notion of graphs on surfaces, namely Hurwitz galaxies (see e.g. \cite{duchi2014bijections} or \cite{PJ}).
\begin{Definition}
 \label{def:hurwitzgalacxy}
 A Hurwitz galaxy of type $(g,\mu,\nu)$ is a graph $G$ on an oriented surface $S$ of genus $g$, such that for $m=2g-2+\ell(\mu)+\ell(\nu)$:
 \begin{enumerate}[(i)]
  \item $S\backslash G$ is a disjoint union of open disks,
  \item $G$ partitions $S$ into $\ell(\mu)+\ell(\nu)$ disjoint faces,
  \item these faces may be coloured black and white, such that $\ell(\nu)$ many faces are coloured black and $\ell(\mu)$ many faces are coloured white, such that each edge is incident to a white face on one side and to a black face on the other side,
  \item the white (resp. black) faces are labeled by $1,\dots,\ell(\mu)$ (resp. $1,\dots,\ell(\nu)$), such that a face labeled $i$ is bounded by $\mu_i\cdot m$ vertices,
  \item the vertices are labeled cyclically counterclockwise with respect to the adjacent white faces by $1,\dots,m$,
  \item for each $i\in\{1,\dots,m\}$, there are $m-1$ vertices labeled $i$, which are 2-valent and one vertex labeled $i$, which is 4-valent.
 \end{enumerate}

 An isomorphism between two Hurwitz galaxies is an orientation-preserving homeomorphism of their respective surfaces, which induces an isomorphism of graphs, that preserves vertex- and face-labels.
\end{Definition}
  \begin{figure}
 \begin{center}
   \scalebox{0.9}{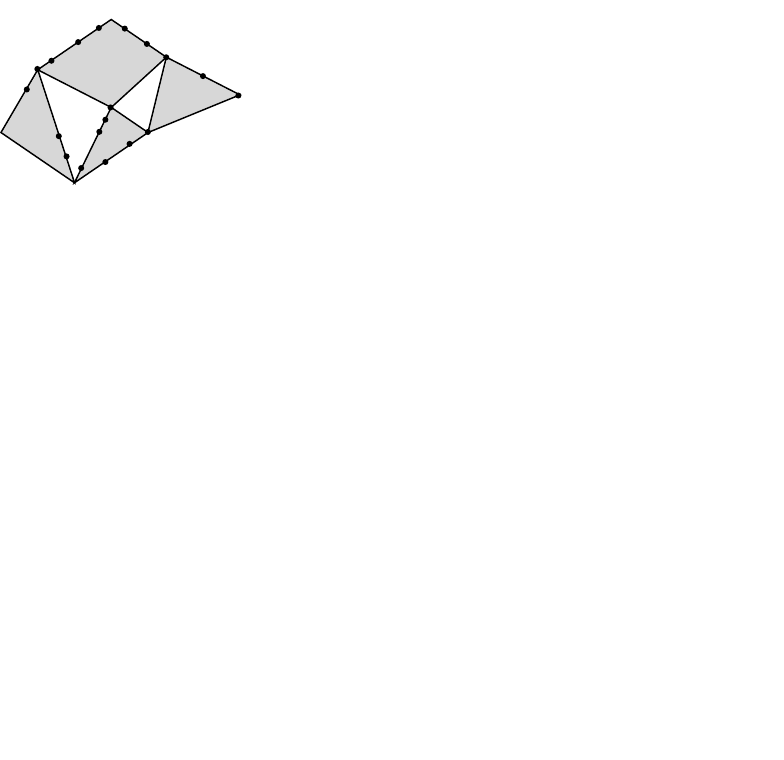}
 \end{center}
 \caption{}
 \label{galaxybranching}
 \end{figure}
 
\begin{proposition}
\label{prop:bij}
 There is a bijection:   \begin{align*}
  \left\{\begin{array}{cl} 
    \textrm{Branching graphs of type }(g,\mu,\nu)\\
  \end{array}\right\}&\leftrightarrow 
  \left\{\begin{array}{cl}  
  \textrm{Hurwitz galaxies of type }(g,\mu,\nu)\\
  \end{array}\right\}.\\
  \end{align*}
\end{proposition}

\begin{proof}
We start with a Hurwitz galaxy $G$ of type $(g,\mu,\nu)$. Draw a vertex in each white face und connect this vertex to the vertices surrounding this face. Now remove the vertices of the old graph $G$. We obtain a branching graph of type $(g,\mu,\nu)$ by distributing the labels naturally. Obviously, we may reverse this process and thus get the bijection as desired.
\end{proof}

\begin{Example}
We illustrate the construction in the proof of Proposition \ref{prop:bij} in Figure \ref{galaxybranching}. We start with a Hurwitz galaxy of type $(0,(2,1,3),(1,2,1,2))$ and obtain the corresponding branching graph of type $(0,(2,1,3),(1,2,1,2))$. The green numbers display the labels of the faces of the galaxy and the labels of the faces and vertices of the branching graph.
\end{Example}

We will construct Hurwitz covers $f$ from branching graphs $\varGamma$. In this construction, we will actually use Hurwitz galaxies $G$. Moreover, we want to relate the automorphism groups. To be more preice, we will see that there are natural bijections between the set of Hurwitz covers of type $(g,\mu,\nu)$, branching graphs of type $(g,\mu,\nu)$ and Hurwitz galaxies of type $(g,\mu,\nu)$. Furthermore, we will see that for a Hurwitz cover $f$, the corresponding branching graph $\varGamma$ and Hurwitz galaxy $G$, there are natural isomorphisms between their automorphism groups.\\
We note that only branching graphs of type $(g,(d),(d))$ have automorphisms. This may be seen by an easy graph theoretic argument. We will give a proof by connecting the automorphisms of branching graphs to automorphisms of factorizations in the symmetric group in Section \ref{sub:poly}.

We can compute Hurwitz numbers in terms of isomorphism classes of branching graphs of type $(g,\mu,\nu)$. We denote the set of all isomorphism classes of branching graphs of type $(g,\mu,\nu)$ by $\mathcal{B}_g(\mu,\nu)$ .

\begin{proposition}[\cite{OPsecondversion},\cite{GJV},\cite{PJ}]
\label{prop:hurwitzgraphcorr}
 With notation as above, we have:
 \[\mathcal{H}_g(\mu,\nu)=\sum_{\varGamma\in \mathcal{B}_g(\mu,\nu)}\frac{1}{|\mathrm{Aut}(\varGamma)|}\]
\end{proposition}

\begin{figure}
\begin{center}
 \scalebox{0.4}{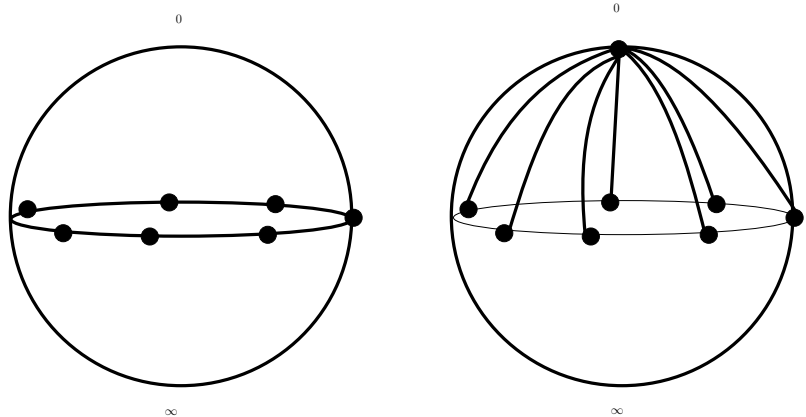}
 \caption{On the left, graph on the sphere, whose pullback yields a Hurwitz galaxy. On the right, graph on the sphere, whose pullback yields a branching graph.}
 \label{corr}
  \end{center}
\end{figure}

The idea behind the proof of Proposition \ref{prop:hurwitzgraphcorr} is to express Hurwitz galaxies and branching graphs as pullbacks of certain graphs on $\PC$ in the following sense: Fix some $f\in\mathbb{H}_g(\mu,\nu)$. Draw the graph whose vertices are the $m=2g-2+\ell(\mu)+\ell(\nu)$ roots of unity and whose edges connect them as in the left graph in Figure \ref{corr}. The pre-image of this graph under $f$ is a Hurwitz galaxy of type $(g,\mu,\nu)$ and each Hurwitz galaxy of type $(g,\mu,\nu)$ appears that way. Similar for branching graphs, we draw the graph whose vertices are the $m$ roots of unity and $0$ on $\PC$ and whose edges connect $0$ to each root of unity as in the right graph in Figure \ref{corr} and take the pre-image.

\section{Pruned double Hurwitz numbers}
\label{chapter:PrunedHurwitz}
In this section, we present our results on pruned double Hurwitz numbers. We begin by defining these objects and formulate our first main result, namely the equivalence between double Hurwitz numbers and pruned double Hurwitz numbers. This theorem expresses double Hurwitz numbers as a weighted sum over pruned double Hurwitz numbers of the same genus. The rest of this section is devoted to proving this theorem.

\label{sub:corr}
As in \cite{DN13} we define the set $\mathcal{PB}_g(\mu,\nu)$ of pruned branching graphs of type $(g,\mu,\nu)$ to be the subset of $\mathcal{B}_g(\mu,\nu)$ consisting of all branching graphs of type $(g,\mu,\nu)$ without leaves. This leads to our main definition, which we introduce here generalizing the definition of pruned simple Hurwitz numbers in \cite{DN13}. 
 \begin{Definition}
 Let $\mu,\nu$ be partitions of the same positive integer $d$. Let $g$ be a non-negative integer. We define the \textit{pruned double Hurwitz number} to be
 \[\mathcal{PH}_g(\mu,\nu):=\sum_{\varGamma\in \mathcal{PB}_g(\mu,\nu)}\frac{1}{|\mathrm{Aut}(\varGamma)|}.\]
 Sometimes, we don't care about automorphisms. Thus we define the \textit{modified pruned double Hurwitz number} to be
 \[\widehat{\mathcal{PH}}_g(\mu,\nu):=\sum_{\varGamma\in\mathcal{PB}_g(\mu,\nu)}1.\]
 \end{Definition}
 \begin{Remark}
  Contrary to the double Hurwitz number $\mathcal{H}_g(\mu,\nu)$, the pruned double Hurwitz number $\mathcal{PH}_g(\mu,\nu)$ is not symmetric in $\mu$ and $\nu$. For example $\mathcal{PH}_0(\mu,(d))=0$, since there are no pruned trees (we exclude the isolated vertex by convention). However, $\mathcal{PH}_0((2),(1,1))=1$.
 \end{Remark}

 By our discussion about automorphisms in Section \ref{chapter:Preliminaries}, we have \[\mathcal{PH}_g(\mu,\nu)=\widehat{\mathcal{PH}}_g(\mu,\nu),\] whenever $0$ or $\infty$ is not fully ramified.\vspace{\baselineskip}\\
In fact we may express the double Hurwitz number as a weighted sum over certain modified pruned double Hurwitz numbers of smaller degree (we have to take the modified Hurwitz numbers, since removing vertices might introduce unwanted automorphisms). The idea is, that we iteratively remove all leaves of the branching graphs until none are left. To make our main result precise, we have to introduce some notation.
\begin{Definition}
 Let $\sigma$ be some rooted forest on the vertex set $\{1,\dots,n\}$. We define $\mathrm{deg}(i)$ to be the number of successors of the vertex $i$. Moreover, we define $\Delta(\sigma)=(\mathrm{deg}(1),\dots,\mathrm{deg}(n))$ to be the ordered degree sequence of $\sigma$. If $\sigma$ has $k$ components, we call $\Delta(\sigma)$ a degree sequence of type $(n,k)$.
\end{Definition}
Note, that some $n$-tuple $(a_1,\dots,a_n)$ is the ordered degree sequence of some rooted forests on $n$ vertices and $k$ components if and only if $\sum a_i=n-k$.
\begin{theorem}
\label{thm:mainresult}
 Let $n=\ell(\nu)$ and let $\mu,\nu$ be partitions of the same positive integer $d$. Then we get:
  \begin{align*}
  \mathcal{H}_g(\mu,\nu)=&\sum_{\tilde{\nu}_1=1,\dots,\tilde{\nu}_n=1}^{\nu_1,\dots,\nu_n}\sum_{\substack{I\subset\{1,\dots,\ell(\mu)\}\}}}\widehat{\mathcal{PH}}_g(\mu_I,\tilde{\nu})\cdot\\&\left(\sum_{\substack{I_1\sqcup\dots\sqcup I_n=I^c:\\\left|\mu_{I_i}\right|=\nu_i-\tilde{\nu}_i}}\dbinom{2g-2+\ell(\mu)+n}{2g-2+\ell(\mu_I)+n,\ell(\mu_{I_1}),\dots,\ell(\mu_{I_n})}\cdot\right.\\
  &\left(\prod_{s=1}^n\ell(\mu_{I_s})!\right)\cdot\left.\sum_{\substack{(\Delta_1,\dots,\Delta_n)\\\Delta_i\textrm{ degree sequence}\\\textrm{of type }(\tilde{\nu}_i+|I_i|,\tilde{\nu}_i) }}\prod_{i=1}^n\right.\\
  &\left.\left(\sum_{j=1}^{\tilde{\nu_i}}\binom{|I_i|-1}{(\Delta_i)_1,\dots,(\Delta_i)_{j-1},(\Delta_i)_j-1,(\Delta_i)_{j+1},\dots,(\Delta_i)_{|\nu_i|+|I_i|}}\right)\cdot\right.\\&\left.\prod_{k\in I_i}(\mu_{I_i})_k^{(\Delta_i)_k}\right)
 \end{align*}
\end{theorem}

\begin{Remark}
We note that if $|\mu_I|\neq\tilde{\nu}$, then $\widehat{\mathcal{PH}}(\mu_I,\tilde{\nu})=0$.\\
Moreover, by inverting the relation we see that pruned Hurwitz numbers are determined by their classical counterparts as well.
\end{Remark}

\begin{Example} Before we start with the proof of Theorem \ref{thm:mainresult}, we give some examples. The Hurwitz numbers appearing in this example were computed with GAP procedures which can be found on \href{https://sites.google.com/site/marvinanashahn/computer-algebra}{https://sites.google.com/site/marvinanashahn/computer-algebra}.
\label{ex:ev1}
 \begin{enumerate}
  \item  We would expect, that \[\mathcal{H}_1((2,2),(1,1,1,1))=\mathcal{PH}_1((2,2),(1,1,1,1))\] and indeed our program yields \[\mathcal{H}_1((2,2),(1,1,1,1))=\mathcal{PH}_1((2,2),(1,1,1,1))=23040.\]
  \item For $\mathcal{H}_0((2,2,2),(3,2,1))$ the formula yields
  \begin{align*}\mathcal{H}_0((2,2,2),(3,2,1))=&\mathcal{PH}_0((2,2,2),(3,2,1))+\\&3\cdot\mathcal{PH}_0((2,2)(2,1,1))\cdot0+\\&3\cdot\mathcal{PH}_0((2,2),(1,2,1))\cdot4.
  \end{align*}
  Indeed our computations show
  \begin{align*}
   &\mathcal{H}_0((2,2,2),(3,2,1))=1728\\
   &\mathcal{PH}_0((2,2,2),(3,2,1))=1152\\
   &\mathcal{PH}_0((2,2),(1,2,1))=48,
  \end{align*}
  which verifies our formula.
 \end{enumerate}
 \end{Example}

Now we may define a construction similar to the construction in the proof of Proposition 3.4 in \cite{DN13}. Firstly, we introduce some new notation: Let $\mu$ be an ordered partition and let $I\subset\{1,\dots,\ell(\mu)\}$, then we denote $\mu_I=(\mu_i)_{i\in I}$.\\
 The following construction associates a pruned branching graph to a branching graph in algorithmic way. We exclude the case $\ell(\nu)=1$, i.e. the case of trees, since in this case our algorithm leaves a single vertex and by convention we excluded this case.
\begin{construction}
\label{constr:prunedfaceforest}
 Let $\varGamma$ be a branching graph of type $(g,\mu,\nu)$, such that $\ell(\nu)>1$. We now construct a subgraph of $\varGamma$ which will indeed be a pruned branching graph of type $(g,\tilde{\mu},\tilde{\nu})$, such that $\tilde{\mu}\subset\mu$, $1\le\tilde{\nu}_i\le\nu_i$ and $|\tilde{\mu}|=|\tilde{\nu}|$.
 \begin{enumerate}
  \item We remove all leaves of $\varGamma$. That is, we remove the vertices of valency $1$, all adjacent half-edges and the adjacent full-edge. Moreover, we remove all half-edges with the same label as the removed full-edge in the whole graph.
    \item After that, we relabel the edges, such that the labels form a set of the form $\{1,..,k\}$ for some $k$.
  \item If the resulting graph $\tilde{\varGamma}$ is pruned, the process stops, if not, we start again.
 \end{enumerate}
 When this process stops, we obtain a pruned branching graph $\tilde{\varGamma}$ of some type $(g,\tilde{\mu},\tilde{\nu})$ with $\tilde{\mu}$ and $\tilde{\nu}$ as above. We call $\tilde{\varGamma}$ the underlying pruned branching graph of $\varGamma$.\\
 Note that we may perform this process for each face seperately. For a face $F$, we call the resulting face $\tilde{F}$ the underlying pruned face.
 \end{construction}

 We refer to Construction \ref{constr:prunedfaceforest} as \textit{pruning}. The resulting underlying pruned branching graph is unique.

\begin{Definition}
\label{def:forestoftype}
 Let $v$ and $\tilde{\nu}$ be integers with $v\geq\tilde{\nu}$ and let $F$ be a rooted forest with $v$ vertices and $\tilde{\nu}$ components. Moreover, let the non-root vertices be bilabeled by some set $I$ and some set $E$, i.e. each non-root vertex has two labels. Let the root-vertices be labeled by some set $R$, such that $v-\tilde{\nu}=|E|=|I|$, $|R|=\tilde{\nu}$. We call $F$ a forest of type $(\tilde{\nu},I,E,R)$.\\
 If we drop the labeling by the set $E$, we call $F$ a forest of type $(\tilde{\nu},I,R)$.
\end{Definition}

\begin{proposition}
\label{prop:bijection}
 Let $\nu$ and $n$ be positive integers and fix some positive integer $m$. Moreover, let $\mathcal{E}$ be some set of order $k$ contained in $\{1,\dots,m\}$. There is a weighted bijection
  \begin{align*}
  \left\{\begin{array}{cl} 
    \textrm{Faces }F\textrm{ of branching}\\
    \textrm{graphs on }n\\
    \textrm{vertices with perimeter }\nu\\
    \textrm{and with full-edge}\\
    \textrm{labels in }\mathcal{E}
  \end{array}\right\}&\leftrightarrow 
  \left\{\begin{array}{cl}  
  \textrm{Triples }(\tilde{F},\varGamma,\mu)\textrm{, such that}\\
  \tilde{F}\textrm{ is a pruned face}\\\textrm{of a branching graph}\\
  \textrm{with perimeter }\tilde{\nu}\leq\nu,\\
  \varGamma\textrm{ is forest of type}\\
  (\tilde{\nu},I,E,R(\tilde{\nu})), \textrm{ for some}\\
  I\subset\{1,\dots,n\},\ E\subset\mathcal{E},\\|
  I|=|E|\textrm{ and}\\
  \textrm{an ordered partition }\mu\textrm{, such that }\\
  \ell(\mu)=|I|
  \textrm{and }\tilde{\nu}+|\mu|=\nu\\
  \end{array}\right\}.\\
  \end{align*}
\end{proposition}

While the proof of this proposition involves some intricate combinatorics, the idea is rather simple: Starting with the face $F$ of the branching graph, we associate a pruned face $\tilde{F}$ as in Construction \ref{constr:prunedfaceforest}. Considering the graph induced by $F-\tilde{F}$, i.e. removing the underlying pruned face, we obtain a forest $\varGamma$. For the other direction, starting with a pruned face $\tilde{F}$ and a forest $\varGamma$, there are several ways of reconstructing a face $F$ by gluing the forst into the pruned face.

\begin{proof}
 We give an algorithm for each direction of the bijection. Let $F$ be a face of a branching graph with a total of $m$ edges, such that $F$ has perimeter $\nu$ with underlying pruned face $\tilde{F}$ of perimeter $\tilde{\nu}$. Furthermore, let $I$ be the set of vertex-labels and $E$ the set of edge-labels not contained in $\tilde{F}$ but in $F$. Let $\mu$ be the partition of the perimeters of those vertices we remove in the pruning process, such that the entries of $\mu$ are labeled by $I$, i.e. the vertex labeled $i$ has perimeter $\mu_i$. We see immediately that $|\mu|=\nu-\tilde{\nu}$, since we remove $|\mu_I|\cdot m$ edges, where we count all full-edges twice, except the ones incident to the underlying pruned face, which we count once. We construct a forest of some type $(\tilde{\nu},I,E,R)$.
 \begin{enumerate}[(a)]
  \item By definition each label occurs exactly $\tilde{\nu}$ times in $\tilde{F}$, such that we can divide the boundary of $\tilde{F}$ in $\tilde{\nu}$ many segments, such that each segment is incident to an edge with a given label exactly once. By convention, each segment starts with the label $1$. We label the segments cyclically counterclockwise by $R_1,\dots,R_{\tilde{\nu}}$, where we assign $R_1$ to the segment containing the full-edge with the smallest label in the face.
  \item Now we contract these segments in $F$ to a root vertex, one for each of the $\tilde{\nu}$ many components. We relabel these components by reassigning each edge label to the adjacent vertex which is further away from the root vertex. This yields the set $E$. The root vertex is labeled by its segment, which corresponds to the set $R$.
  \item Moreover we contract all half-edges.
 \end{enumerate}
 Furthermore, each non-root vertex is by definition labeled by $I$, thus we obtain a forest of type $(\tilde{\nu},I,E,R)$ as above. This construction is unique. \vspace{\baselineskip}\\
 For the other direction, we start with a tuple $(\tilde{F},\varGamma)$, such that $\tilde{F}$ has perimeter $\nu$ and $\varGamma$ is a forest of type $(\tilde{\nu},I,E,R)$. We start by labeling the segments of the boundary of $\tilde{F}$ as above by $R_1,\dots,R_{\tilde{\nu}}$ cyclically counterclockwise, such that the segment labeled $R_1$ contains the full-edge with the smallest label. Now,  we glue the forest into the pruned face as follows:
 \begin{enumerate}
  \item We give the forest $\varGamma$ the canonical orientation. We label each edge by the label of its target-vertex corresponding to the set $E$.
  \item We introduce a partial ordering on the edges of $\varGamma$ in the following way: For two edges $e,e'$ we define $e\ge e'$, if they are contained in the same tree and $e$ lies on the unique path from the respective root vertex to $e'$.
  \item Relabel the edges of $\tilde{F}$ by $\mathcal{E}-E$, such that the edge labeled $i$ is relabeled by the $i-th$ element of $\mathcal{E}-E$ in the natural order. Add the half-edges in $\{1,\dots ,m\}-(\mathcal{E}-E)$ to the pruned face, such that half-edges are labeled cyclically counterclockwise by $\{1,\dots,m\}$.
  \item Attach the maximal edges adjacent to the vertex labeled $R_i$ to the segment labeled $R_i$ as follows: Let $e$ be an edge of the forest adjacent to the root vertex labeled $R_i$ with target vertex labeled by $(h,j)\in I\times E$. Glue an edge labeled $j$ to the half-edge labeled $j$ in the segment $R_i$, label the new vertex of valency $1$ by $h$ and add $\mu_h\cdot m-1$ half-edges to $h$ that are cyclically labeled by $\{1,\dots,m\}$. Thus, each edge label occurs $\mu_h$ times at $h$. This procedure is unique.
  \item Remove the maximal elements from the ordering.
  \item Attach the maximal edges in the new ordering as follows: Let $e$ be such an edge of $\tilde{F}$, such that the source vertex of $e$ is labeled by $(h_s,j_s)\in I\times E$ and the target vertex by $(h_t,j_t)\in I\times E$. Glue an edge labeled $j_t$ to a half-edge labeled by $j_t$ that is adjacent to the vertex labeled $h_s$ in the new face. There are $\mu_{h_s}$ many half-edges labeled $j_t$ adajacent to $h_s$, thus we have $\mu_{h_s}$ many choices for this edge and thus $\mathrm{per}(v)^{\mathrm{val}(v)-1}$ choices for each vertex. Label the new vertex of valency $1$ by $h_t$ and add $\mu_{h_t}\cdot m-1$ many half-edges to $h_t$ as in step 4.
  \item Repeat steps $5$ and $6$ iteratively for the other edges and vertices of the forest.
 \end{enumerate}
 We obtain a face $F$ of perimeter $\nu$. One can check that both constructions are inverse to each other. The choices in step $6$ are the only choices we have and thus we obtain a weighted bijection as desired.
\end{proof}

\begin{Example}
We use the construction in Proposition \ref{prop:bijection} in the example Figure \ref{forestconstr}. We start with a pruned face with perimeter $12$. We remove vertices with labels $5-11$ and edges with labels $2,6-11$. The remaining labels $1,2,3,5$ are relabeled as $1,2,3,4$. We obtain a pruned face with perimeter $4$, the rooted forest in Figure \ref{forestconstr} and the partition $(1,1,2,1,1,1,1)$. These objects satisfy all conditions.
 \begin{figure}
 \begin{center}
   \scalebox{0.75}{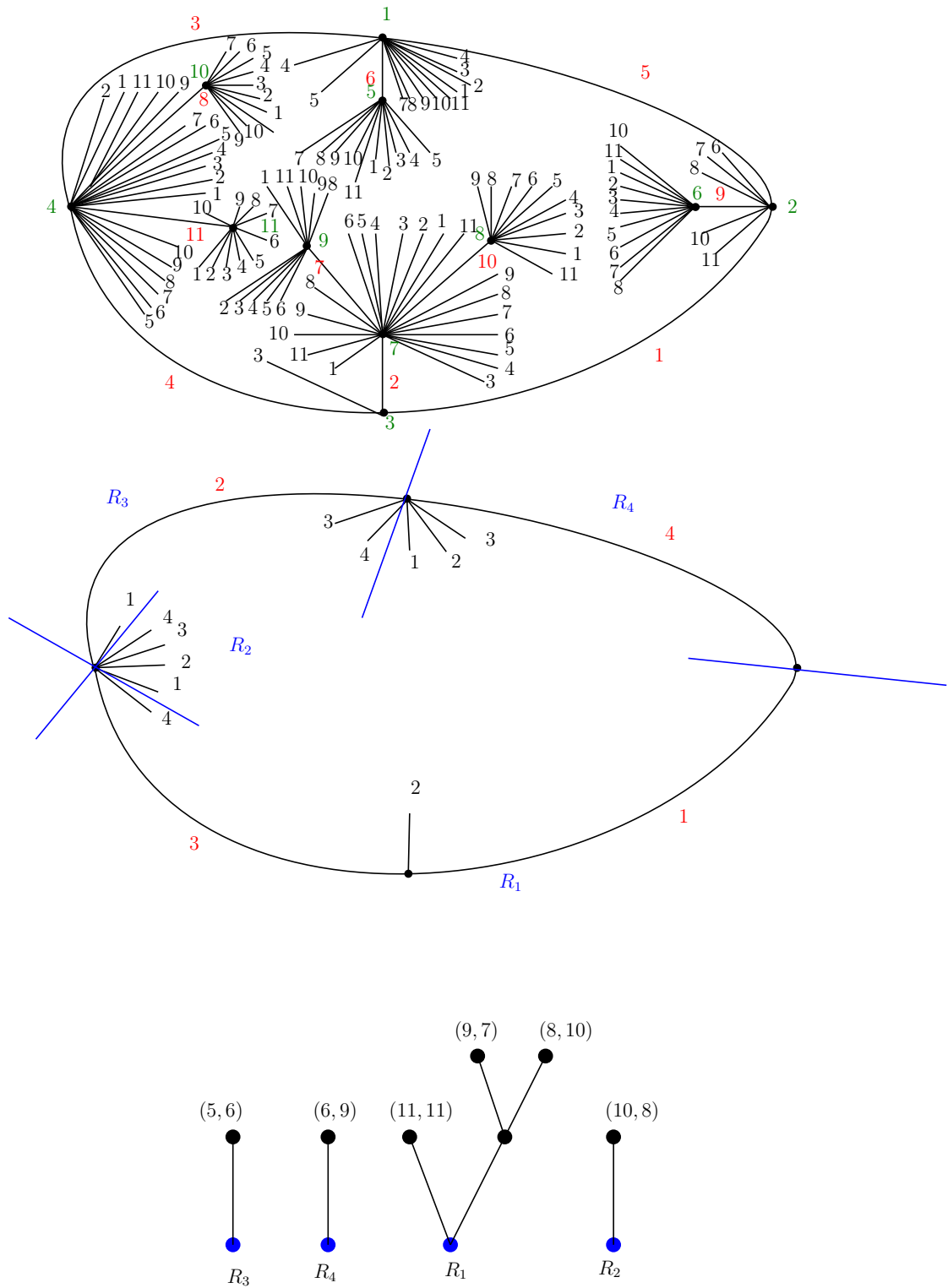}
 \end{center}
 \caption{The red labels correspond to the full-edges and the green labels to the vertices.}
 \label{forestconstr}
\end{figure}
\end{Example}

\begin{proposition}
\label{thm:classicprunedcorr}
 Let $n=\ell(\nu)$ and let $\mu,\nu$ be partitions of the same positive integer $d$. Then we get:
 \begin{align*}
  \mathcal{H}_g(\mu,\nu)=	&\sum_{\tilde{\nu}_1=1,\dots,\tilde{\nu}_n=1}^{\nu_1,\dots,\nu_n}\sum_{\substack{I\subset\{1,\dots,\ell(\mu)\}\textrm{,}\\ \textrm{such that}\\ \left|\mu_I\right|=\left|\tilde{\nu}\right|}}\widehat{\mathcal{PH}}_g(\mu_I,\tilde{\nu})\cdot\\&\left(\sum_{\substack{I_1\sqcup\dots\sqcup I_n=I^c:\\\left|\mu_{I_i}\right|=\nu_i-\tilde{\nu}_i}}\dbinom{2g-2+\ell(\mu)+\ell(\nu)}{2g-2+\ell(\tilde{\mu})+\ell(\tilde{\nu}),\ell(\mu_{I_1}),\dots,\ell(\mu_{I_n})}\cdot\right.\\&\left(\prod_{i=1}^n\ell(\mu_{I_i})!\right)\cdot
  \left.\sum_{\substack{(\varGamma_1,\dots,\varGamma_n):\\ \varGamma_i \textrm{ is a rooted forest of type}\\ (\tilde{\nu_i},I_i,R(\tilde{\nu}_i)))}}\prod_{k=1}^n\prod_{\substack{v \textrm{ non-root}\\ \textrm{vertex of }\varGamma_k}}(\mu_{I_k})_v^{\mathrm{val}(v)}\right),
 \end{align*}
 where $R(\tilde{\nu}_i)$ is the index set $R(\tilde{\nu}_i)=\{R_1,\dots,R_{\tilde{\nu}_i}\}$.
\end{proposition}

\begin{proof}
 The proof is similar to the proof of Proposition 3.4 in \cite{DN13}. The given formula is a weighted sum over pruned branching graphs. As already seen in Construction \ref{constr:prunedfaceforest}, we may assign a unique pruned branching graph to each branching graph. For the other direction we apply Proposition \ref{prop:bijection} to each face iteratively. Recall that we may obtain a branching graph of type $(g,\mu,\nu)$ from a pruned branching graph of type $(g,\mu_I,\tilde{\nu})$ for some $I$, such that $1\le\tilde{\nu}_i\le\nu_i$ and $|\mu_I|=|\tilde{\nu}|$. We can do this by choosing a decomposition $I^c=I_1\sqcup\dots\sqcup I_{\ell(\nu)}$, such that $|\mu_{I_i}|=\nu_i-\tilde{\nu}_i$ and adding vertices to the face labeled $i$, whose perimeters correspond to $\mu_{I_i}$, in a tree-like manner. Thus, adding $\ell(\mu_{I_i})$ vertices means adding just as many edges.\\
 The desired formula may be reformulated as follows: There is weighted bijection
   \begin{align*}
  \left\{\begin{array}{cl} 
    \textrm{Branching graph }\varGamma\\
    \textrm{of type }(g,\mu,\nu)\\
  \end{array}\right\}&\leftrightarrow 
  \left\{\begin{array}{cl}  
  \textrm{Tuple }(\varGamma,I,(I_1,\dots,I_n),(\varGamma_1,\dots,\varGamma_n)),\\
  \textrm{such that}\varGamma\textrm{ is a pruned branching graph}\\
  \textrm{of type }(g,\mu_I,\tilde{\nu})\textrm{ for some subset }I,\\
  I_1\sqcup\dots\sqcup I_n=I^c\textrm{, such that}\\
  \varGamma_i \textrm{ is a rooted forest of type}\\ (\tilde{\nu_i},I_i,R(\tilde{\nu}_i))\\
  \textrm{and }\tilde{\nu}_i-\nu_i=|\mu_{I_i}|\\
  \end{array}\right\}.\\
  \end{align*}
 Now we count the number of branching graphs of type $(g,\mu,\nu)$ with underlying pruned branching graph of type $(g,\tilde{\mu},\tilde{\nu})$. We do this by reconstructing branching graphs of type $(g,\mu,\nu)$ from pruned branching graphs of type $(g,\mu_I,\tilde{\nu})$:
 
 \vspace{\baselineskip}
 Fix a pruned branching graph $\varGamma$ of type $(g,\mu_I,\tilde{\nu})$ for some $I\subset\{1,\dots,\ell(\mu)\}$, such that $|\mu_I|=|\tilde{\nu}|$. We need to add vertices and edges as described above. Firstly, we distribute the perimeters of the vertices to the faces, that means, we choose some decomposition $I^c=I_1\sqcup\dots\sqcup I_n$, such that $|\mu_{I_i}|=\nu_i-\tilde{\nu_i}$. Moreover, we distribute the edge-labels of the pruned branching graph as well as the set of edge labels, we add to face $i$, i.e. we choose a decomposition of the $2g-2+\ell(\mu)+\ell(\nu)$ edge labels $E(\varGamma)=\tilde{E}\sqcup E_1\sqcup\dots\sqcup E_n$, such that $|\tilde{E}|=2g-2+\ell(\mu_I)+\ell(\tilde{\nu})$ and $|E_i|=|I_i|$.\\
 Now we may add vertices and edges as described to construct some branching graph of type $(g,\mu,\nu)$. For each branching graph constructed that way, the face $i$ contracts to some forest of type $(\tilde{\nu}_i,I_i,E_i,R(\tilde{\nu}_i))$ as in Proposition \ref{prop:bijection}. As noted in Proposition \ref{prop:bijection}, each forest of type $(\tilde{\nu}_i,I_i,E_i,R(\tilde{\nu_i}))$ corresponds to \[\prod_{\substack{v \textrm{ non-root}\\ \textrm{vertex of }\varGamma_k}}(\mu_{I_i})_v^{\mathrm{val}(v)-1}\] many different faces. Thus we obtain
  \begin{align*}
  \mathcal{H}_g(\mu,\nu)=&\sum_{\tilde{\nu}_1=1,\dots,\tilde{\nu}_n=1}^{\nu_1,\dots,\nu_n}\sum_{\substack{I\subset\{1,\dots,\ell(\mu)\}\textrm{,}\\ \textrm{such that}\\ \left|\mu_I\right|=\left|\tilde{\nu}\right|}}\widehat{\mathcal{PH}}_g(\mu_I,\tilde{\nu})\cdot
  \sum_{\substack{I_1\sqcup\dots\sqcup I_n=I^c:\\\left|\mu_{I_i}\right|=\nu_i-\tilde{\nu}_i}}\sum_{\substack{\{1,\dots,m\}=\\ \tilde{E}\sqcup E_1\sqcup\dots\sqcup E_n \textrm{,}\\\textrm{such that}\\|\tilde{E}|=2g-2+\ell(\mu_I)+\ell(\tilde{\nu})|\\\textrm{and}\\|E_i|=\ell(\mu_{I_i})}}\\ &\sum_{\substack{(\varGamma_1,\dots,\varGamma_n):\\ \varGamma_i \textrm{ is a rooted forest of type}\\ (\nu_i,\tilde{\nu_i},I_i,E_i,R(\tilde{\nu}_i)))}}\prod_{k=1}^n\prod_{\substack{v \textrm{ non-root}\\ \textrm{vertex of }\varGamma_k}}(\mu_{I_i})_v^{\mathrm{val}(v)-1}.
 \end{align*}
 To obtain the formula we want to prove, it is enough to observe, that a different choice of edge labels does not change the factor\[\prod_{k=1}^n\prod_{\substack{v \textrm{ non-root}\\ \textrm{vertex of }\varGamma_k}}(\mu_{I_i})_v^{\mathrm{val}(v)-1}\] in the formula above and that there are \[\dbinom{2g-2+\ell(\mu)+\ell(\nu)}{2g-2+\ell(\mu_I)+\ell(\tilde{\nu}),\ell(\mu_{I_1}),\dots,\ell(\mu_{I_n})}\] ways to choose the edge labels on the underlying pruned branching graph, as well as the set of $\ell(\mu_{I_i})$ edge labels to add to the face $i$. Moreover, there are \[\prod_{s=1}^n\ell(\mu_{I_s})!\] ways to distribute the edge labels to the vertices of each graph. Thus we obtain the desired formula.
\end{proof}

This is a generalization of the respective theorem in \cite{DN13} in the sense, that for double Hurwitz numbers we obtain a weighted count over tuples of forests, where in the simple Hurwitz numbers case, each tuple is counted with weight $1$.
In fact, we may simplify the formula in Theorem \ref{thm:classicprunedcorr}, by using the following result on the number of rooted forests.

\begin{theorem}
 Let $S\subset\{1,\dots,n\}$ be a fixed set and let $\mathcal{T}_{n,S}$ be the set of rooted forests with $n$ vertices and $|S|$ components, such that the roots are labeled by $S$.
 \[\sum_{\sigma\in\mathcal{T}_{n,S}}x_1^{\mathrm{deg}(1)}\cdots x_n^{\mathrm{deg}(n)}=(x_1+\dots+x_n)^{n-|S|-1}\sum_{i\in S}x_i\]
\end{theorem}

\begin{proof}
 See for example \cite{RS} page 29.
\end{proof}

Thus for some fixed degree sequence $\Delta=(\delta_1,\dots,\delta_n)$ of type $(n,k)$, the number of rooted forests on $k$ components, such that the roots are labeled by some set $S\subset\{1,\dots,n\}$, is:
\begin{align*}&\textrm{Coefficient of }x_1^{\delta_1}\cdots x_n^{\delta_n}\textrm{ in} \sum_{\sigma\in\mathcal{T}_{n,S}}x_1^{\mathrm{deg}(1)}\cdots x_n^{\mathrm{deg}(n)}=\\&\sum_{i\in S}\binom{n-|S|-1}{\delta_1,\dots,\delta_{i-1},\delta_i-1,\delta_{i+1},\dots,\delta_n}\end{align*}
Using this result, we see that for a fixed partition $\mu$ and for each degree sequence $(\delta_1,\dots,\delta_n)$ of type $(\ell(\mu)+|S|,|S|)$, there are \[\sum_{i\in S}\binom{\ell(\mu)+|S|-|S|-1}{\delta_1,\dots,\delta_{i-1},\delta_i-1,\delta_{i+1},\dots,\delta_n}=\sum_{i\in S}\binom{\ell(\mu)-1}{\delta_1,\dots,\delta_{i-1},\delta_i-1,\delta_{i+1},\dots,\delta_n}\]
many forests of type $(n-|\mu|,\{1,\dots,\ell(\mu)\},R(n-|\mu|)$, that correspond to the factor \[\prod_{\substack{v \textrm{ non-root}\\ \textrm{vertex of }\varGamma_k}}(\mu_{I_i})_v^{\mathrm{val}(v)-1}=\prod_{i=1}^{\ell(\mu)}\mu_i^{\mathrm{\delta_i}}\] in the formula in Proposition \ref{thm:classicprunedcorr}.
Thus, by adjusting the set $\{1,\dots,n\}$ to $I_i\sqcup R(\tilde{\nu}_i)$ in the formula and choosing $S=R(\nu_i)$ (without loss of generality $R(\nu_i)$ corresponds to the first $\tilde{\nu_i}$ entries of each degree sequence), we obtain Theorem \ref{thm:mainresult}.

\section{A cut-and-join recursion for pruned Hurwitz numbers}
\label{sub:rec}
In the standard case double Hurwitz numbers admit a recursion. We now formulate and prove an analogue for the pruned case. This is a generalization of Proposition 3.2 in \cite{DN13}, where a similar cut-and-join recursion is proved for pruned simple Hurwitz numbers.

\begin{theorem}
\label{thm:prunedrecursion}
  Let $\mu$ and $\nu$ be partitions of the same positive integer $d$ and $g$ a non-negative integer, such that $\ell(\mu)+\ell(\nu)+2g-2>0$ and $(g,\ell(\nu))\neq(0,1)$ and $(g,\ell(\nu))\neq(0,2)$. Then the following recursion formula holds
  \begin{alignat*}{4}
  \mathcal{PH}_g(\mu,\nu)&=\frac{1}{2}\sum_{i=1}^{\ell(\nu)}\sum_{I\subset\{1,\dots,\ell(\mu)\}}\sum_{\substack{\alpha+\beta+|\mu_{I^c}|\\=\nu_i}}&&(|I^c|+1)\cdot\frac{(m-1)!}{(m-(|I^c|+1))!}\cdot\prod_{a\in\mu_{I^c}}a\cdot\\& &&\alpha\cdot\beta\cdot\widehat{\mathcal{PH}}_{g-1}(\mu_I,(\nu_{S\backslash\{i\}},\alpha,\beta))\\
  &+\frac{1}{2}\sum_{i=1}^{\ell(\nu)}\sum^{stable}_{\substack{g_1+g_2=g,}}\sum_{\substack{J_1\sqcup J_2=\\S\backslash\{i\}}}\sum_{\substack{I_1,I_2\subset\\\{1,\dots,\ell(\mu)\}\\\textrm{disjoint}}} \sum_{\substack{\alpha+\beta+\\|\mu_{(I_1\sqcup I_2)^c}|\\=\nu_i}}&&\frac{(m-1)!}{m_1!m_2!}\cdot\alpha\cdot\beta\cdot\\
  & &&(|(I_1\sqcup I_2)^c|+1)!\cdot\prod_{a\in\mu_{(I_1\sqcup I_2)^c}}a\cdot\\
  & &&\widehat{\mathcal{PH}}_{g_1}(\mu_{I_1},(\nu_{J_1},\alpha))\cdot\widehat{\mathcal{PH}}_{g_2}(\mu_{I_2},(\nu_{J_2},\beta))\\
  &+\sum_{i<j}\sum_{I\subset\{1,\dots,\ell(\mu)\}}\sum_{\substack{\alpha+|\mu_{I^c}|=\\\nu_i+\nu_j}}&&\frac{(m-1)!}{(m-(|I^c|+1))!}\cdot\prod_{a\in\mu_{I^c}}a\cdot\alpha\cdot(|I^c|+1)!\cdot\\& &&\widehat{\mathcal{PH}}_g(\mu_I,(\nu_{S\backslash\{i,j\}},\alpha)),
  \end{alignat*}
  where $m_1=2g_1-2+|I_1|+|J_1|$, $m_2=2g_2-2+|I_2|+|J_2|$, $S=\{1,\dots,\ell(\nu)\}$ and the term \textit{stable} in the sum expresses that we exclude terms with $(g_i,|J_i|)=(0,2)$ or $(g_i,|J_i|)=(0,1)$.
\end{theorem}
  \begin{figure}
 \begin{center}
  \scalebox{0.5}{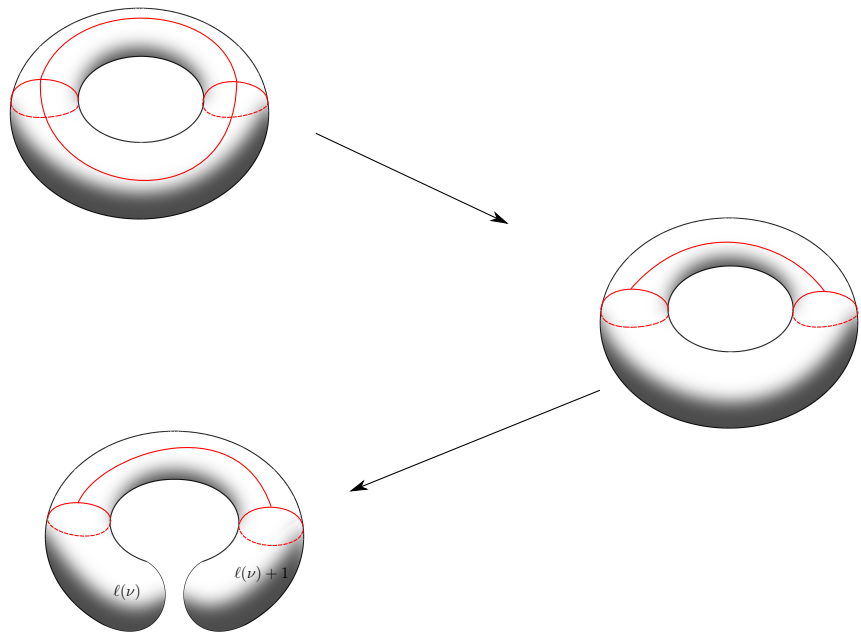}
 \end{center}
 \caption{We start with a graph of some type $(g,\mu,\nu)$ (drawn red). After removing $m$, we lose genus. The degeneration process is visualized.}
 \label{firstcase}
\end{figure}

\begin{Remark}
For each term on the right hand side of the equation, the number of simple ramification points is smaller than $2g-2+\ell(\mu)+\ell(\nu)$. Thus, the base cases for this recursion are $\widehat{\mathcal{PH}}_0(\mu,\nu)$, where $\ell(\mu),\ell(\nu)\le2$.
\end{Remark}

 \begin{Example}
Before we begin by computing some examples verifying our formula. (We implemented a procedure in GAP \cite{GAP4} computing the numbers below, which can be found in \href{https://sites.google.com/site/marvinanashahn/computer-algebra}{https://sites.google.com/site/marvinanashahn/computer-algebra}.)
 \label{ex:ev2}
\begin{enumerate}
 \item Our first example for the recursion is for $\mathcal{PH}_0((3),(1,1,1))$. The formula yields: \[\mathcal{PH}_0((3),(1,1,1))=3\cdot\mathcal{PH}_0((3),(1,2))\cdot 2.\] In fact, the computations show, that 
 \begin{align*}
 &\mathcal{PH}_0((3),(1,1,1))=6\\
 &\mathcal{PH}_0((3),(1,2))=1.
 \end{align*}
 \item Our second example is for $\mathcal{PH}_1((2,2),(1,1,2))$. The recursion yields: \begin{align*}\mathcal{PH}_1((2,2),(1,1,2))=&\mathcal{PH}_0((2,2),(1,1,1,1))\cdot\frac{1}{2}+\mathcal{PH}_1((2,2),(2,2))\cdot2+\\&2\cdot\mathcal{PH}_1((2,2),(1,3))\cdot3+4\cdot\widehat{\mathcal{PH}}_1((2),(1,1))\cdot16.\end{align*}
 Our computations yield \begin{align*}
  &\mathcal{PH}_1((2,2),(1,1,2))=1920\\ 
  &\mathcal{PH}_0((2,2),(1,1,1,1))=576\\
  &\mathcal{PH}_1((2,2),(2,2))=160\\
  &\mathcal{PH}_1((2,2),(1,3))=1248\\
  &\mathcal{PH}_1((2),(1,1))=1.                                                                                                                                                                                                                                                                   \end{align*}
These numbers satisfy the recursion as expected.
\end{enumerate}
 \end{Example}

\begin{figure}
 \begin{center}
  \scalebox{0.6}{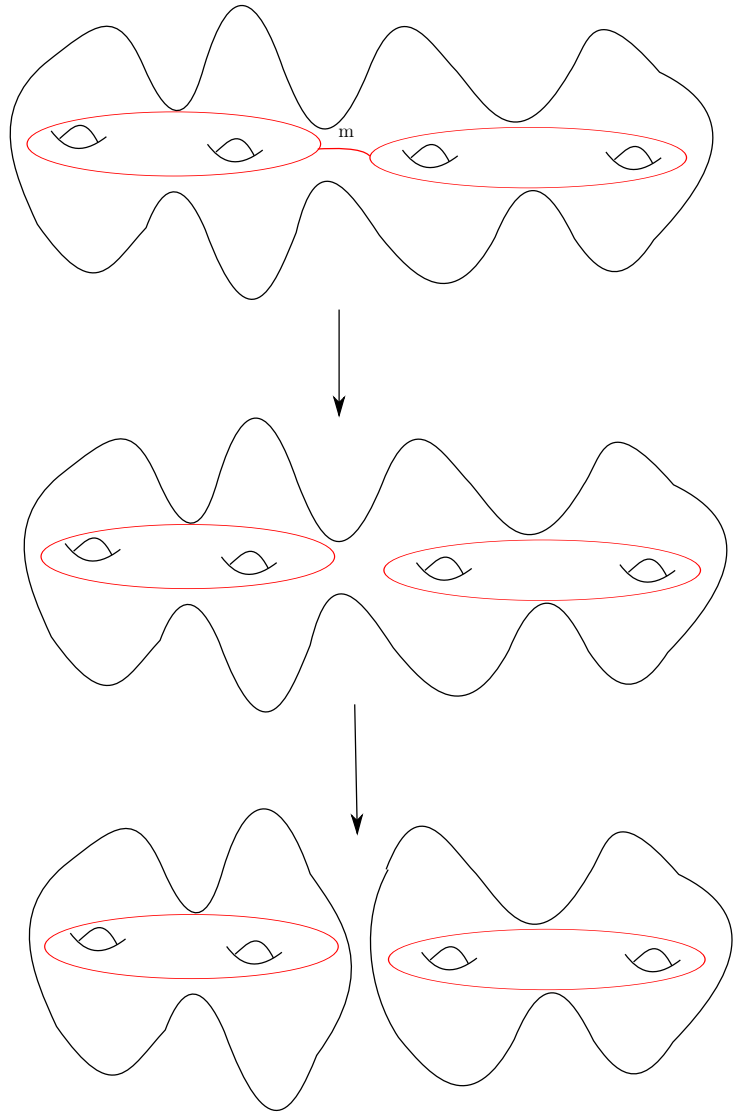}
 \end{center}
 \caption{We start with a graph of some type $(g,\mu,\nu)$. After removing $m$, the graph decomposes in two graphs and we degenerate.}
 \label{secondcase}
\end{figure}

\begin{figure}
 \begin{center}
 \scalebox{0.45}{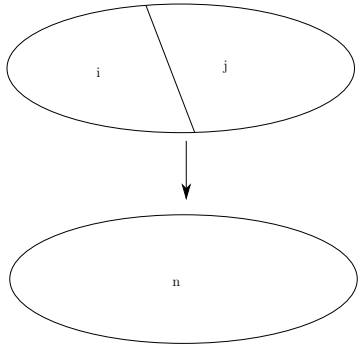}
 \end{center}
 \caption{We start with a graph of some type $(g,\mu,\nu)$. After removing $m$, two faces join to a new one.}
 \label{thirdcase}
 \end{figure}

\begin{center}
\begin{figure}
\scalebox{0.7}{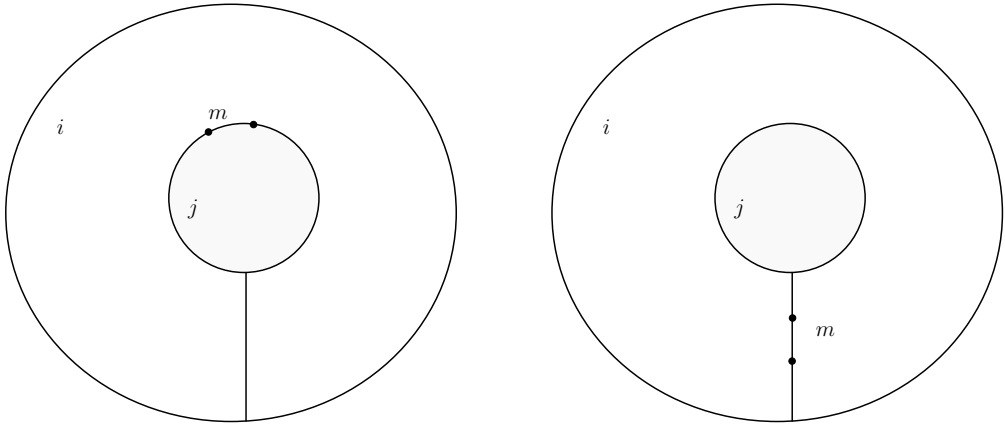}
\caption{Removing the edge labeled $m$ in the left picture corresponds to the third case in the proof of Theorem \ref{thm:prunedrecursion}. However, reconstructing as in the third case, allows placing the edge labeled $m$ as in the right picture, which actually corresponds to the second case.}
\label{fig:stablerecursion}
\end{figure}
\end{center}

The idea behind this recursion is similar to the one in \cite{DN13}, which we aim to generalize. We start with a branching graph of type $(g,\mu,\nu)$ and remove the full-edge labeled $m$ and all half-edges with the same label. This may leave a graph that is not pruned. In that case, we apply Construction \ref{constr:prunedfaceforest} and obtain a new pruned graph $\varGamma$. We exclude the cases, where $\ell(\nu)\leq2$, since our procedure is not well-defined in the case, where the graph we start with is just a cycle. Since the graph is pruned, the removed edges either form a path or look locally like the left graph in Figure \ref{fig:stablerecursion}. We can classify the possible cases for the new graph:
\begin{enumerate}
\item The new branching graph obtained that way is a pruned branching graph of type $(g-1,\mu_I,(\nu_{S\bs\{i\}},\alpha,\beta))$ for some subset $I\subset\{1,\dots,\ell(\mu)\}$, $i\in\{1,\dots,\ell(\nu)\}$ and $\alpha,\beta>0$, such that $\alpha+\beta+|\mu_{I^c}|=\nu_i$. Note, that we require for $\varGamma$ in order to be a branching graph, that its faces are homeomorphic to open disks. Thus, we need to degenerate the surface, $\varGamma$ is embedded on, as illustrated in Figure \ref{firstcase}.
\item The new graph is a disjoint union of two pruned branching graphs of respective type $(g_1,\mu_{I_1},(\nu_{J_1},\alpha))$ and $(g_2,\mu_{I_2},(\nu_{J_2},\beta))$, whereas $J_1\sqcup J_2=S\backslash\{i\}$ and $I_1,I_2\subset\{1,\dots,\ell(\mu)\}$, such that $g_1+g_2=g$, $I_1\cap I_2=\emptyset$ and $\alpha+\beta+|\mu_{(I_1\sqcup I_2)^c}|\\=\nu_i$. This is illustrated in Figure \ref{secondcase}.
\item The new graph is a pruned branching graph of type $(g,\mu_I,(\nu_{S\bs\{i,j\}},\alpha))$, where $i,j\in\{1,\dots,\ell(\nu)\}$ and $I\subset\{1,\dots,\ell(\mu)$, such that $\alpha+|\mu_{I^c}|=\nu_i+\nu_j$. This is illustrated in Figure \ref{thirdcase}.
\end{enumerate}
Now, we give algorithms to reconstruct graphs of type $(g,\mu,\nu)$ from graphs in each of the three cases.

\begin{algorithm}
\label{alg:1}
 We begin this algorithm by fixing $\varGamma$ to be some pruned branching graph of type $(g-1,\mu_I,(\nu_{S\bs\{i\}},\alpha,\beta))$ as in the first case. First we need to embedd $\varGamma$ on a surface of genus $g$, such that the faces labeled $\ell(\nu)$ and $\ell(\nu)+1$ are joined, reversing the second step in Figure \ref{firstcase}. We construct a pruned branching graph of type $(g,\mu,\nu)$ as follows, reversing the first step in Figure \ref{firstcase}:
 \begin{enumerate}
  \item Set $T=\{1,\dots,m\}$, $U=I^c$.
  \item Choose an edge label $k$ in $T$ and attach an edge with that label to the face labeled $\ell(\nu)$ of perimeter $\alpha$.
  \item Set $T\rightarrow T\bs\{k\}$.
  \item Choose a vertex label $l$ in $U$ and attach a vertex of perimeter $\mu_l$ to the other end of the edge, we attached in step $2$.
  \item Choose an edge label $k$ in $T$ and attach an edge with that label to the vertex, we just attached.
  \item Set $T\rightarrow T\bs\{k\}$ and $U\rightarrow U\bs\{l\}$.
  \item Repeat steps $3-5$ until $|U|=\emptyset$.
  \item Attach the last edge we attached to the path to the face labeled $\ell(\nu)+1$ of perimeter $\beta$.
  \item Relabel the edges of the graph without the new path by $T$, such that the order of the edge labels is maintained.
  \item Label the face obtained by joining $\ell(\nu)$ and $\ell(\nu)+1$ by $i$ and adjust the labels of the other faces.
  \end{enumerate}
  The new graph obtained that way is a pruned branching graph of type $(g,\mu,\nu)$.
\end{algorithm}

\begin{algorithm}
\label{alg:2}
 We begin by fixing $\varGamma_1$ and $\varGamma_2$ to be some pruned branching graphs of respective type $(g_1,\mu_{I_1},(\nu_{J_1},\alpha))$ and $(g_2,\mu_{I_2},(\nu_{J_2},\beta))$ as in the second case. First, we need to embedd those graph of a surface of genus $g$, such that the face labeled $|J_1|+1$ of $\varGamma_1$ and the face labeled $|J_2|+1$ of $\varGamma_2$ are joined, reversing the second step in Figure \ref{secondcase}. We construct a pruned branching graph of type $(g,\mu,\nu)$ as follows, reversing the first step in Figure \ref{secondcase}:
 \begin{enumerate}
 \item Set $T=\{1,\dots,m\}$, $U=(I_1\sqcup I_2)^c$.
 \item Choose an edge label $k$ in $T$ and attach an edge with that label to the face labeled $|J_1|+1$ of $\varGamma_1$ of perimeter $\alpha$.
  \item Set $T\rightarrow T\bs\{k\}$.
  \item Choose a vertex label $l$ in $U$ and attach a vertex of perimeter $\mu_l$ to the other end of the edge, we attached in step $2$.
  \item Choose an edge label $k$ in $T$ and attach an edge with that label to the vertex, we just attached.
  \item Set $T\rightarrow T\bs\{k\}$ and $U\rightarrow U\bs\{l\}$.
  \item Repeat steps $3-5$ until $|U|=\emptyset$.
  \item Attach the last edge we attached to the path to the face labeled $|J_2|+1$ of $\varGamma_2$ of perimeter $\beta$, joining the two graphs.
  \item Relabel the edges of the graph without the new path by $T$, such that the order of the edge labels is maintained.
  \item Label the new face obtained by joining both graphs by $i$ and adjust the labels of the other faces.
  \end{enumerate}
  The new graph obtained that way is a pruned branching graph of type $(g,\mu,\nu)$.
\end{algorithm}

\begin{algorithm}
\label{alg:thirdcase}
 We begin by fixing $\varGamma$ to be some pruned branching graph of type \[(g,\mu_I,(\nu_{S\bs\{i,j\}},\alpha))\] as in the third case. We construct a pruned branching graph of type $(g,\mu,\nu)$ as follows, reversing the process in Figure \ref{thirdcase}.
  \begin{enumerate}
 \item Set $T=\{1,\dots,m\}$, $U=I^c$.
 \item Choose an edge label $k$ in $T$ and attach an edge with that label to the face labeled $\ell(\nu)-1$ of $\varGamma$ of perimeter $\alpha$.
  \item Set $T\rightarrow T\bs\{k\}$.
  \item Choose a vertex label $l$ in $U$ and attach a vertex of perimeter $\mu_l$ to the other end of the edge, we attached in step $2$.
  \item Choose an edge label $k$ in $T$ and attach an edge with that label to the vertex, we just attached.
  \item Set $T\rightarrow T\bs\{k\}$ and $U\rightarrow U\bs\{l\}$.
  \item Repeat steps $3-5$ until $|U|=\emptyset$.
  \item Attach the last edge we attached to the path to the new face labeled $\ell(\nu)-1$ in such a way, that it is divided in two faces of respective perimeter $\nu_i$ and $\nu_j$.
  \item Relabel the edges of the graph without the new path by $T$, such that the order of the edge labels is maintained.
  \item If $\nu_i\neq\nu_j$ label the face of perimeter $\nu_i$ by $i$ and the face of perimter $\nu_j$ by $j$. If $\nu_i=\nu_j$ choose a way to label the faces by $i$ and $j$. Adjust the labels of the other faces.
  \end{enumerate}
  The new graph obtained that way is a pruned branching graph of type $(g,\mu,\nu)$.
\end{algorithm}

In all three algorithms we have to make some choices, thus the result of each algorithm is not uniquely determined by the initial conditions. The next step in order to prove Theorem \ref{thm:prunedrecursion} is to analyze the number of choices we have in each algorithm. However, in each algorithm not every resulting graph will yield the pruned branching graph we began with, after removing the edge labeled $m$ and pruning. In the first two algorithms the graphs, where $m$ lies on the path we added will fulfil this property. In the third algorithm, we allowed the path to join itself in the last step. Thus allowing $m$ on the whole path isn't enough as illustrated in Figure \ref{fig:stablerecursion}. However, we will repair this below in the proof of Theorem \ref{thm:prunedrecursion}.\vspace{\baselineskip}\\
We call the resulting graphs with the edge labeled $m$ on the path we attached, the \textit{relevant graphs}.

\begin{lemma}
\label{lem:recursion}
The number of relevant graphs of type $(g,\mu,\nu)$, we obtain
 \begin{enumerate}
  \item in Algorithm \ref{alg:1} from a fixed graph of type $(g-1,\mu_I,(\nu_{S\bs\{i\}},\alpha,\beta))$ is\[\alpha\cdot(|I^c|+1)!\cdot\frac{(m-1)!}{(m-(|I^c|+1))!}\cdot\prod_{a\in\mu_{I^c}}a\cdot\beta,\]
  \item in Algorithm \ref{alg:2} from two fixed graphs with each of respective type $(g_1,\mu_{I_1},(\nu_{J_1},\alpha))$ and $(g_2,\mu_{I_2},(\nu_{J_2},\beta))$ is 
  \item in Algorithm \ref{alg:thirdcase} from a fixed graph of type $(g,\mu_I,(\nu_{S\bs\{i,j\}},\alpha))$ is \[\alpha\cdot(|I^c|+1)!\cdot\frac{(m-1)!}{(m-(|I^c|+1))!}\cdot\prod_{a\in\mu_{I^c}}a.\]
 \end{enumerate}

\end{lemma}

\begin{proof}
 \begin{enumerate}
  \item In the first case, there are $\alpha$ many ways to attach the first edge. There are $|I^c|!$ many ways to distribute the vertex labels to the path. Moreover, since we only count relevant graphs, we have $|I^c|+1$ possibilities to assign the label $m$ to some edge on the path. After assigning the label $m$, there are $m-1$ many labels to assign to the $m-(|I^c|+1)-1$ edges on the path without a label, which yields a factor of $\frac{(m-1)!}{(m-(|I^c|+1))!}$. When we attach an edge to a vertex label $i\in I^c$, there are $\mu_i$ many ways to attach that edge in each step. Thus we obtain a factor of $\prod_{a\in\mu_{I^c}}a$. Finally, no graph occurs twice in this construction, thus we proved the first statement.
  \item The second case works analogously to the first one.
  \item In the third case, the factors occur the same way as in the first and second case, except for the eighth and tenth step in Algorithm \ref{alg:thirdcase}. If $\nu_i\neq\nu_j$ in the eighth step, we have two choices to attach the last edge to the face and only one possibility in the tenth step. If $\nu_i=\nu_j$, we have only one choice in the eighth step, but two choicesin the tenth step. This would yield a factor of $2$. However, the algorithm produces each graph twice by the following argument: If the path is not attached to itself, we cannot distinguish which end of the graph was attached to the face first. If the path is attached to itself, one vertex of the path is trivalent and two adjacent edges are contained in a cycle. We cannot distinguish which of those two edges was attached last. This yields a factor of $\frac{1}{2}$ and the third statement is proved.
  \end{enumerate}
\end{proof}

Now, we are ready to finish the proof of Theorem \ref{thm:prunedrecursion}.

\begin{proof}[Proof of Theorem \ref{thm:prunedrecursion}]
The three reconstructive algorithms produce all graphs of type $(g,\mu,\nu)$. We need to make sure, that each graph is obtained only once. However, we have already seen, that is not true, since the third algorithm produces graphs that contribute to the second case, as illustrated in Figure \ref{fig:stablerecursion}. However, those are exactly the graphs of the second case, where one graph is of type $(0,\tilde{\mu},\tilde{\nu})$, such that $\ell(\tilde{\nu})=2$. Thus, we just exclude those cases in the second algorithm. We can also exclude those graphs with $(0,\ell(\tilde{\nu}))=(0,1)$ since $\widehat{\mathcal{PH}}_0(\mu_I,\tilde{\nu})=0$.\vspace{\baselineskip}\\
Moreover, if $\alpha=\beta$ in the first case, we may switch the labeling of the respective faces and the first algorithm yields the same relevant graphs. Thus, we have to adjust the count by $\frac{1}{2}$, if $\alpha=\beta$. However, for $\alpha\neq\beta$ the first algorithm yields the same relevant graphs for graphs in $\mathcal{PH}_{g-1}(\mu_I,(\nu_{S\backslash\{i\}},\alpha,\beta))$ as in $\mathcal{PH}_{g-1}(\mu_I,(\nu_{S\backslash\{i\}},\beta,\alpha))$, since the construction is symmetric in $\alpha$ and $\beta$. Thus, we adjust those summands by a factor of $\frac{1}{2}$ as well.\vspace{\baselineskip}\\
A similar argument accounts for the factor $\frac{1}{2}$ in the second case and the recursion follows.
\end{proof}

\section{Polynomiality of pruned double Hurwitz numbers and connection to the symmetric group}
\label{sub:poly}
It is well known, that double Hurwitz numbers in arbitrary genus are piecewise polynomial in the $\mu_i$ and $\nu_i$. The first proof was given in \cite{GJV}. The proof for pruned double Hurwitz numbers works analogously. We start by recalling the structure of the proof in \cite{GJV}: We fix some tuple $(g,\mu,\nu)$. There are only finitely many branching graphs of that type. In each branching graph $\varGamma$ of type $(g,\mu,\nu)$ we drop the half-edges and obtain a new graph $\tilde{\varGamma}$, which we call the \textit{skeleton} of $\varGamma$. For each type $(g,\mu,\nu)$, there are only finitely many skeletons, which may be obtained from such a branching graph. However, many branching graphs may have the same skeleton. We define $S(g,\mu,\nu,\tilde{\varGamma})$ to be the number of branching graphs of type $(g,\mu,\nu)$ with skeleton $\tilde{\varGamma}$. Thus, we may compute $\mathcal{H}_g(\mu,\nu)$ as weighted sum over all skeletons, where each skeleton $\tilde{\varGamma}$ is weighted by $S(g,\mu,\nu,\tilde{\varGamma})$. This is a finite sum, since all but finitely many skeletons will be weighted by $0$. In \cite{GJV} it was proved that $S(g,\mu,\nu,\tilde{\varGamma})$ behaves piecewise polynomially in the entries of $\mu$ and $\nu$ by using Erhart theory and that each polynomial has degree $4g-3+\ell(\mu)+\ell(\nu)$. Thus by refining the hyperplanes, piecewise polynomiality follows for $\mathcal{H}_g(\mu,\nu)$.\\
This approach is feasible for pruned double Hurwitz numbers, since the property of a branching graph being pruned is inherent in its skeleton. Thus, $\mathcal{PH}_g(\mu,\nu)$ may be computed as a weighted sum over all pruned skeletons, where each skeleton $\tilde{\varGamma}$ is weighted by $S(g,\mu,\nu,\tilde{\varGamma})$. The piecewise polynomiality follows analogously. The precise statement is as follows:

\begin{theorem}
\label{thm:poly}
 Let $k$ and $l$ be two positive integers and $g$ some non-negative integer, then $\mathcal{PH}_g(\mu,\nu)$ is piecewise polynomial in the entries of $\mu$ and $\nu$ (where $\ell(\mu)=k$ and $\ell(\nu)=l$) of degrees up to $4g-3+k+l$. The ``leading'' term of degree $4g-3+k+l$ is non-zero, i.e. $\mathcal{PH}_g(t\mu,t\nu)$, as a function in $t$ of positive integer values, is a polynomial of degree $4g-3+k+l$.
\end{theorem}

In order to make the contributions of each skeleton more precise we introduce the notion of a reduced branching graph, which will also make the results concerning the connection to the symmetric group easier.
\begin{Definition}
 For a branching graph $\varGamma$ let $\varGamma^s$ be the graph obtained from $\varGamma$ by dropping all its half-edges. We call $\varGamma^s$ the skeleton of $\varGamma$. 
\end{Definition}

\begin{Notation}
 Let $\varGamma$ be an edge-labeled graph on a surface. We define a corner of the skeleton to be a tuple $(v,e,e',F)$, such that $e$ and $e'$ are both full-edges adjacent to $v$ and $F$ and $e'$ is positioned after $e$ counterclockwise.\\
 We call a corner \textit{descendant}, if the label of $e'$ is smaller than the label of $e$.
\end{Notation}

\begin{Definition}
 Let $d$ and $g$ be positive integers, moreover let $\mu$ and $\nu$ be ordered partitions of $d$. We define a \textit{reduced branching graph of type }$(g,\mu,\nu)$ to be a graph $\varGamma$ on an oriented surface $S$ of genus g, such that for $m=\ell(\mu)+\ell(\nu)-2+2g$:
 \begin{enumerate}[(i)]
 \item $S\backslash\varGamma$ is a disjoint union of open disks,
  \item there are $\ell(\mu)$ vertices, labeled $1,\dots,\ell(\mu)$, such that the vertex labeled $i$, is incident to $\mu_i$ half-edges labeled $p$,
  \item there are exactly $m$ full edges labeled by $1,\dots,m$,
  \item the $\ell(\nu)$ faces are labeled by $1,\dots,\ell(\nu)$ and the face labeled $i$ has perimeter $per(i)=\nu_i$, by which we mean, it contains $\mu_i$ many half-edges labeled $p$,
  \item there is at least one half-edge labeled $p$ in each descending corner.
 \end{enumerate}
Note that we allow loops at the vertices.
\end{Definition}

\begin{Remark}
There is a natural bijection between branching graphs of type $(g,\mu,\nu)$ and reduced branching graphs of type $(g,\mu,\nu)$ given by pulling back an additional edge in the star graph adjacent to $0$ and an unramified point $p$ and forgetting the all half-edges not labeled not labeled $p$ on the source-surface.
\end{Remark}

The contribution of each skeleton is the number of possibilities to distribute half-edges labeled $p$ to each vertex to obtain a reduced monodromy graph, such that the perimeter of the vertex labeled $i$ is $\mu_i$ and the perimeter of the face labeled $j$ is $\nu_j$. We compute the standard and pruned polynomials in one example.

\begin{figure}[H]
\begin{centering}

\scalebox{0.5}{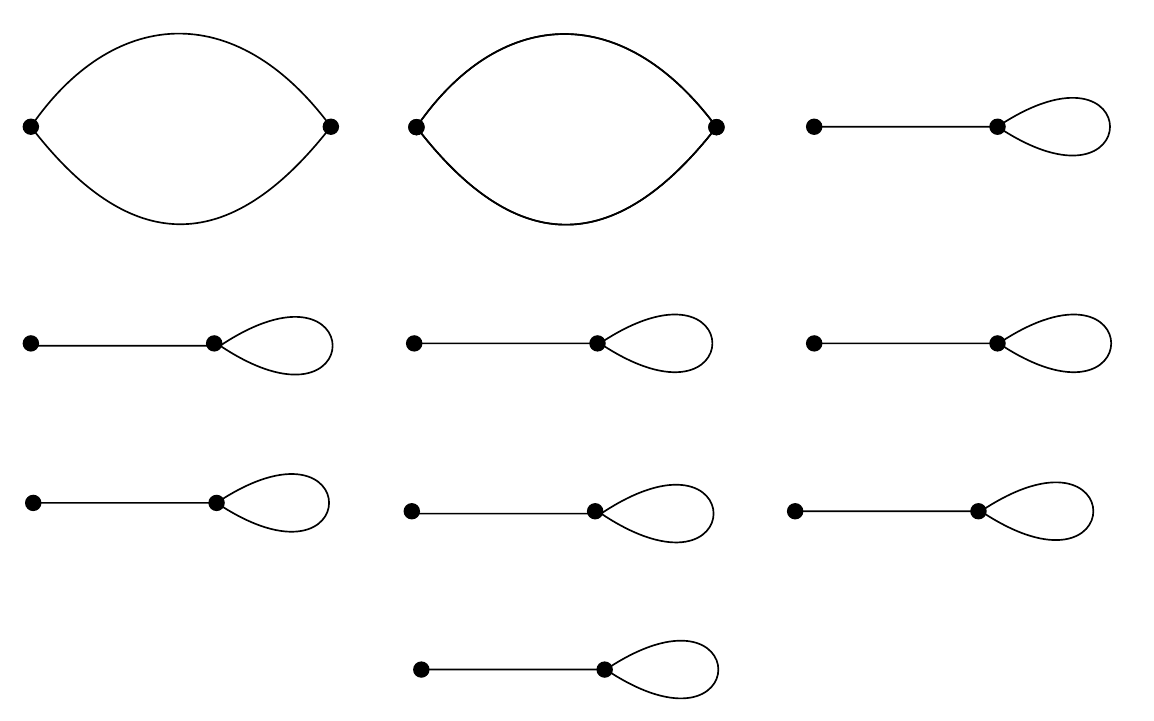}
\caption{The black labels are vertex labels, the green ones are face labels and the red ones are edge labels.}
\label{fig:expoly}
\end{centering}

\end{figure}

\begin{Example}
  We compute the polynomials in genus $0$ for the double Hurwitz numbers $\mathcal{H}_0((a,b),(c,d))$ and their pruned counterparts $\mathcal{PH}_0((a,b),(c,d))$. In this simple case, we can read the contribution directly from the graph without using the procedure of the proof. All possible skeletons are illustrated in Figure \ref{fig:expoly} (in what follows, we enumerate the graphs from the top left to bottom right along the rows). Only the first two are pruned. We compute the polynomial for the chamber $c<a,b<d$.
 \begin{enumerate}
  \item The first two skeletons each contribute a factor of $c$: We need to attach $a$ half-edges labeled $p$ to the vertex labeled $1$ and $b$ half-edges labeled $p$ to the vertex labeled $2$, such that the face labeled $1$ has perimeter $c$ and the face labeled $2$ has perimeter $d$. Since $c<a$, for any $n\in\{1,\dots,c\}$, we can attach $n$ half-edges labeled $p$ to the vertex labeled $1$, such that these half-edges are contained in the face labeled $1$. This determines the entire graph, thus we have $|\{1,\dots,c\}|$ choices.
  \item The third and fourth graph each contribute a factor of $b-c$.
  \item The fifth and sixth graph contribute a factor of $0$.
  \item The seventh and eighth graph each contribute a factor of $a-c$.
  \item The ninth and tenth graph each contribute a factor of $0$.
  Thus we obtain for $c<a,b<d$ 
  \begin{align*}
   &\mathcal{H}_0((a,b),(c,d))=2\cdot c+2\cdot(b-c)+2\cdot(a-c)=2d\\
   &\mathcal{PH}_0((a,b),(c,d))=2\cdot c
  \end{align*}
 \end{enumerate}
 Analogously for the other chambers $d<a,b<c$, $a<c,d<b$ and $b<c,d<a$.

\end{Example}

In Section \ref{chapter:Preliminaries} we have explained the connection between Hurwitz numbers and branching graphs and the connection between Hurwitz numbers and factorizations in the symmetric group. The proof of Theorem \ref{thm:symmetricgroup} yields the following algorithm, which yields the connection between branching graphs and factorizations in the symmetric group. In \cite{PJ}, a similar algorithm is given, which for a given Hurwitz galaxy yields a representation in the symmetric group. However, this algorithm produces the products of permutations $\tilde{\sigma}_i=\tau_i\dots\tau_1\sigma_1$ from which we can recursively deduce $(\sigma_1,\tau_1,\dots,\tau_r,\sigma_2)$. Our algorithm produces $(\sigma_1,\tau_1,\dots,\tau_r,\sigma_2)$ as in Theorem \ref{thm:symmetricgroup} directly and is a direct consequence of the mondromy representation of a branched holomorphic covering.

\begin{Definition}
 Let $\varGamma$ be a reduced branching graph of type $(g,\mu,\nu)$ we call the conjugacy class of the tuple $(\sigma_1,\tau_1,\dots,\tau_m,\sigma_2)$, that is produced by the algorithm below, the \textit{monodromy representation of} $\varGamma$.
\end{Definition}

The notion of a \textit{monodromy representation of a branched covering} in the literature is closely related to the notion defined above. Namely, one can think of a monodromy representation of a cover as a choice of a tuple $(\sigma_1,\tau_1,\dots,\tau_m,\sigma_2)$ as in Theorem \ref{thm:symmetricgroup}. To be more precise: Let $B=\{p_1,p_2,q_1,\dots,q_m\}$ be the set of branch points on $\mathbb{P}^1(\mathbb{C})$, then the tuple $(\sigma_1,\tau_1,\dots,\tau_m,\sigma_2)$ defines a group homomorphism $\Phi:\pi_1(\mathbb{P}^1(\mathbb{C}\backslash B)\to\mathcal{S}_d$. We will see in Proposition \ref{prop:rep} that the monodromy representations of a branched covering and of its corresponding branching graph coincide.

\begin{algorithm}
\label{alg:sym}
 Let $\varGamma$ be a reduced branching graph of type $(g,\mu,\nu)$.
 \begin{enumerate}
  \item Enumerate the half-edges adjacent to vertex $i$ cyclically counterclockwise by \[\sum_{j< i}\mu_j,\dots,\sum_{j\leq i}\mu_j.\] This yields the permutation \[\sigma_1=(1,\dots,\mu_1)\cdots(\sum_{j<\ell(\mu)}\mu_j+1,\dots,\sum_{j\leq\ell(\mu)}\mu_j),\]
  \item Label the $i$-th cycle by $i$,
  \item For the edge labeled $i$, define $\tau_i=(b\ d)$ as in Figure \ref{fig:sym},
  \item Define $(\sigma_2)^{-1}$ to be the permutation whose $i$-th cycle is given by the cyclic numbering of labels of half-edges in the $i$-th face and label the $i$-th cycle by $i$.
  \end{enumerate}
 This gives a tuple $(\sigma_1,\tau_1,\dots,\tau_m,\sigma_2)$ as in Theorem \ref{thm:symmetricgroup}.
\end{algorithm}
\begin{figure}[H]
 \begin{center}
  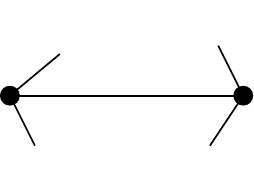
 \end{center}
 \caption{}
 \label{fig:sym}
\end{figure}
Note, that we have a choice in the first step of Algorithm \ref{alg:sym}, namely we didn't specify  where the enumeration starts. However, this just corresponds to conjugations of the resulting monodromy representation, thus the resulting conjugacy class of the algorithm is well-defined.
\begin{Example}
We illustrate Algorithm \ref{alg:sym} for the graph in Figure \ref{examsym}. Note, that we dropped the labels of the vertices and faces for the sake of simplicity. The algorithm yields for the first permutation \[\sigma_1=(1\ 2\ 3\ 4\ 5)(6\ 7)(8\ 9)(10\ 11).\]
The transpositions are\[\tau_1=(8\ 10)\textrm{, }\tau_2=(4\ 7)\textrm{, }\tau_3=(2\ 11)\textrm{ and }\tau_4=(6\ 9)\]
and for the second permutation, we obtain \[(\sigma_2)^{-1}=(1\ 11\ 8\ 6\ 4\ 5)(7\ 9\ 10\ 2\ 3).\]
Indeed, we obtain \[(\sigma_2)^{-1}=\tau_4\tau_3\tau_2\tau_1\sigma_1.\]
 \begin{figure}[H]
  \begin{center}
   \scalebox{0.491}{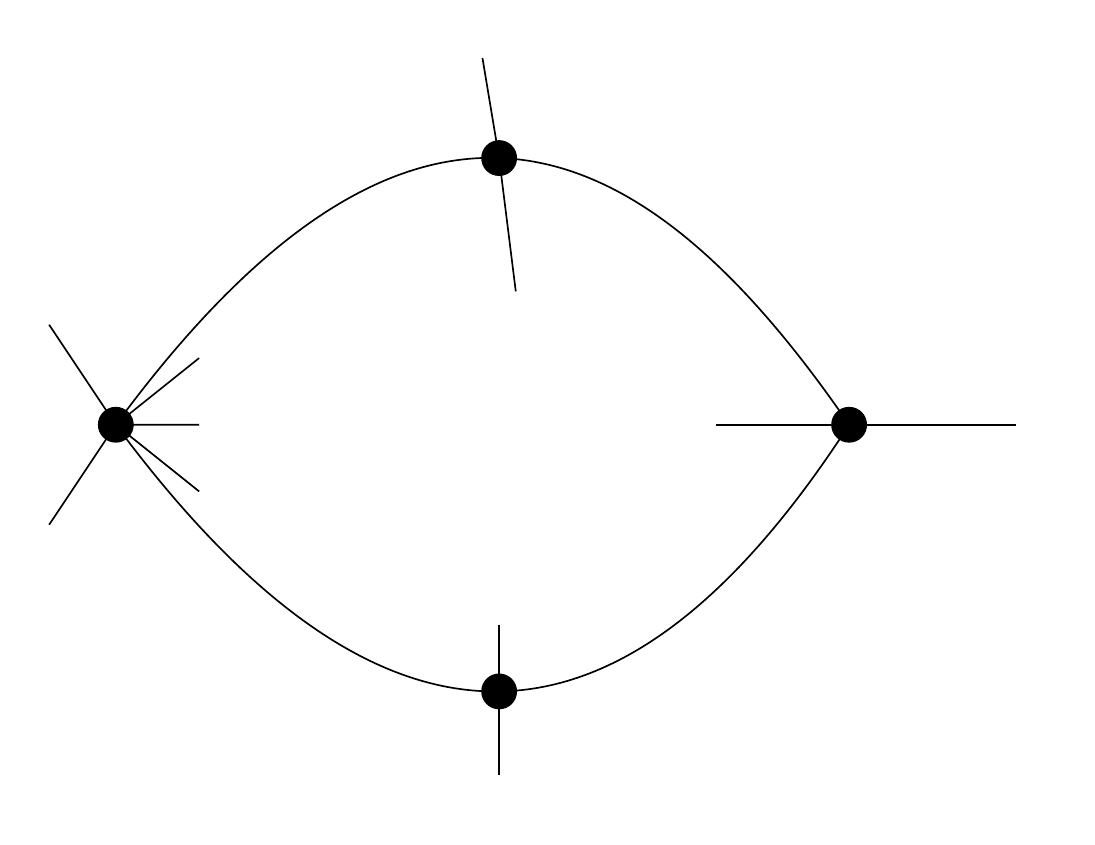}
  \end{center}
  \caption{}
  \label{examsym}
 \end{figure}
\end{Example}

\begin{proposition}
\label{prop:rep}
The monodromy representation of a branched covering and of its corresponding branching graph coincide.
\end{proposition}

The proof is similar to the discussion in Section 4 in \cite{PJ}. Now, we will pick up our discussion about automorphisms in Section \ref{chapter:Preliminaries}. One can check, that two branching graphs $\varGamma$ and $\varGamma'$ are isomorphic, if their corresponding monodromy representations coincide. On the other hand, the conjugation of a tuple in the monodromy representation yields another isomorphic branching graph by relabeling. That means isomorphims between branching graphs correspond to conjugations of the results of Algorithm \ref{alg:sym}. It follows, that automorphisms correspond to conjugations that preserve the result of Algorithm \ref{alg:sym}. However, due to transitvity and the fact that we labeled the disjoint cycles of $\sigma_1$ and $\sigma_2$, it follows that only tuples, where $\sigma_1$ and $\sigma_2$ are $d$-cycles may be invariant under non-trivial conjugations.

We finish this section by giving a classification of pruned Hurwitz numbers in terms of factorizations in the symmetric group, which is an immediate consequence from Algorithm \ref{alg:sym}.

\begin{theorem}
\label{thm:prunsym}
 Let $\mu$ and $\nu$ be partitions of some positive integer $d$. Moreover, let $g$ be some non-negative integer and $m>1$. The following equation holds:
  \begin{alignat*}{3}
   &\mathcal{PH}_g(\mu,\nu)=\\&\frac{1}{d!}\left|
   \begin{cases}
   \begin{rcases}
   &\left(\sigma_1,\tau_1,\dots,\tau_m,\sigma_2\right)\textrm{, such that:}\\ 
   &\bullet\ \sigma_1,\ \sigma_2,\ \tau_i\in\mathcal{S}_d,\\
   &\bullet\ \sigma_2\cdot\tau_m\cdot\dots\cdot\tau_1\cdot\sigma_1=\mathrm{id},\\
   &\bullet\ \mathcal{C}(\sigma_1)=\mu,\ \mathcal{C}(\sigma_2)=\nu\textrm{ and }\mathcal{C}(\tau_i)=(2,1,\dots,1),\\
   &\bullet\ \textrm{The group generated by }\left(\sigma_1,\tau_1,\dots,\tau_m,\sigma_2\right)\\
   &\textrm{acts transitively on }\{1,\dots,d\}\textrm{ and}\\
   &\bullet\ \textrm{The disjoint cycles of }\sigma_1\textrm{ and }\sigma_2\textrm{ are labeled}\\
   &\bullet\ \textrm{For all cycles }\sigma_1^i\textrm{ in }\sigma_1\textrm{there are at least two}\\
   &\textrm{transpositions }\tau_j,\tau_k\textrm{, such that}\\
   &\mathrm{supp}(\tau_j)\cap\mathrm{supp}(\sigma_1^i)\neq\emptyset\neq\mathrm{supp}(\tau_k)\cap\mathrm{supp}(\sigma_1^i),
   \end{rcases}
   \end{cases}
\right|.
  \end{alignat*}
  where for a permutation $\sigma$, we define $\mathrm{supp}(\sigma)$ to be the set of all elements in $\{1,\dots,d\}$, that are not fixed by $\sigma$.
\end{theorem}

\begin{proof}
To begin with, we prove that for each pruned branching graph of type $(g,\mu,\nu)$, Algorithm \ref{alg:sym} produces such a representation. The only condition to check is the last one, but this is immediate, because each cycle $\sigma_1^i$ corresponds to a vertex $i$. This vertex $i$ is not a leaf, because the branching graph we began with is pruned. Thus, there are two edge $e$ and $e'$ adjacent to $i$. However, these edges correspond to two transpositions $\tau_e$ and $\tau_{e'}$, that by construction fulfil the last condition.\\
The other direction follows similarly from the fact, that the monodromy representation of a branching graph is the same as the monodromy representation of the corresponding cover.\\
We excluded $m=1$, due to the fact, that we assume the graph consisting of only one loop and one vertex to be pruned.
\end{proof}

\section*{Acknowledgements}
I owe special thanks to my advisor Hannah Markwig for many helpful conversations, as this paper grew out of the preparation of my Master's thesis. Moreover, I want to thank Norman Do for interesting discussion about his paper \cite{DN13} and an anonymous referee for many helpful remarks on an earlier version of this paper. Thanks are also due to Frank-Olaf Schreyer and Roland Speicher. Finally, I would like to thank Malte Grajewski for his help with Figure \ref{firstcase}.  The author gratefully acknowledges partial support by DFG SFB-TRR 195 "Symbolic tool in mathematics and their applications", project A 14 "Random matrices and Hurwitz numbers" (INST 248/238-1).

\end{document}